\documentclass[12pt]{article}
\usepackage{amsmath}
\usepackage{amssymb}

\usepackage{dsfont}
\usepackage{cite}
\newsymbol \blackbox 1004
\newcommand{\eh}{\hfill}\newlength{\sperr}

\newenvironment{proof}{{\settowidth{\sperr}{\bf\rm
Proof}%
\par\addvspace{0.3cm}\noindent\parbox[t]{1.3\sperr}
{\bf\rm P\eh r\eh o\eh o\eh f\eh }%
}}{\nopagebreak\mbox{}
$\blackbox$\par\addvspace{0.3cm}}

\def\nn{\nonumber}

\def\a{\alpha}
\def\b{\beta}

\def\G{\Gamma}

\def\la{\lambda}
\def\om{\omega}
\def\Om{\Omega}
\def\th{\theta}

\def\vp{\varphi}

\def\wh{\widehat}
\def\wt{\widetilde}
\def\ov{\overline}

\def\BC{{\mathbb C}}

\def\cla{{\mathcal A}}

\def\clm{{\mathcal M}}

\def\clq{{\mathcal Q}}

\def\clt{\mathcal{T}}

\def\rk{{\mathrm{rank}}}

\def\diag{\mathrm{diag}}
\def\col{\mathrm{col}}

\newcommand{\I}{\mathrm{i}}
\newcommand{\1}{\mathbf{1}}

\newtheorem{Pa}{Paper}[section]
\newtheorem{Tm}[Pa]{{\bf Theorem}}
\newtheorem{La}[Pa]{{\bf Lemma}}

\newtheorem{Rk}[Pa]{{\bf Remark}}

\newtheorem{Nn}[Pa]{{\bf Notation}}
\newtheorem{Pn}[Pa]{{\bf Proposition}}

\title{Inversion of the Toeplitz-block Toeplitz matrices and the structure of the corresponding inverse matrices}

\author{Alexander Sakhnovich}

\date{}
\begin{document}
\maketitle

\begin{abstract} 
The results on the inversion of convolution operators and Toeplitz matrices
in the 1-D (one dimensional) case are classical and have numerous
applications. We consider a 2-D case of Toeplitz-block Toeplitz matrices,
describe a minimal information, which is necessary to recover the inverse
matrices, and give a complete characterisation of the inverse matrices.
\end{abstract}




\section{Introduction}\label{Intro}
\setcounter{equation}{0}
The well-known Toeplitz matrices are diagonal-constant matrices whereas
Toeplitz-block Toeplitz matrices $T$ are block Toeplitz matrices, the blocks of which  also are Toeplitz:
\begin{align}\label{T1}&
T=\{\clt_{i-k}\}_{i,k=1}^n, \quad \clt_r=\{t_r^{(j-\ell)}\}_{j,\ell=1}^m.
 \end{align}
Later we use for  Toeplitz-block Toeplitz  matrices the acronym TBT. One may consider TBT-matrices as $2$-D
(two-dimensional) analog of Toeplitz matrices.
 
 Toeplitz matrices (i.e., $1$-D Toeplitz matrices) and their continuous analogs (so called {\it convolution operators} or {\it operators
 with difference kernels}) are very important in mathematical analysis and various applications (see e.g. \cite{BoG, BoeS, Go, Nob, SaSaR, SaL15}
 and references therein).
 The inversion of these matrices and operators is connected with the names of N. Wiener, E. Hopf, N. Levinson, M.G.~Krein,
 I.C. Gohberg, V.A. Ambartsumyan,  L.A. Sakhnovich and many other mathematicians and applied scientists.
 The inversion of convolution operators on a semi-axis (of Wiener-Hopf operators) was studied in numerous papers including 
 the brilliant works \cite{GoKr, Kre}. The situation with the inversion of finite Toeplitz matrices and of convolution operators 
 on a finite interval is (in many respects) more complicated and essentially different from the case of semi-axis.
 Important results on this topic were derived, in particular, in \cite{Amb, GoS, Levin, Shi, Tr}.
 Then, the procedure of recovery of  the operator $T^{-1}$, which is inverse  to a general-type convolution operator $T$ on interval,
 from the action of $T^{-1}$ on two functions  was given in \cite{SaL73} (see also \cite{SaL80, SaL15} and references therein),
 using the method of operator identities.  The same method was applied to the general-type finite Toeplitz matrices in \cite{SaA73}. The structure of the matrices and
 operators $T^{-1}$ was derived in this way as well. The method of operator identities, which was introduced in \cite{SaL68, SaL73},  may be successfully used for the
 inversion of various other structured matrices and operators. See, for instance, \cite{FKRS} and  \cite[Appendix~D]{SaSaR}.
Note also that relations of the form $S^*TS=T$ were used for the study of Toeplitz operators in the seminal work \cite{BrH}.

The inversion  of Toeplitz and Toeplitz-block Toeplitz (TBT) matrices and their continuous analogs is actively studied  in the recent years
as well (see e.g. \cite{Ald, AGKLS, Bin, BogoB, CrL, FKRS, JeVa, NiB, SaL15, Xi} and references therein).
However, in spite of some interesting recent and older works \cite{GeK, Gr, JeVa, Jus, Kalo, Lev, Nap, WaK} on the inversion of TBT-matrices and  of convolution operators in multidimensional spaces, the structure of the corresponding inverse matrices and operators (and the way to recover these inverses from some {\it minimal} information) remained unknown. Our paper deals with this important problem for TBT-matrices
(some less complete results on the inversion of convolution operators on a rectangular are presented in \cite{SaA2017}).

It is easy to show (see Section \ref{Str}) that a TBT-matrix $T$ of the form \eqref{T1} satisfies two matrix identities:
 \begin{align}\label{T2}&
A_pT-TA_p^*=Q_p \quad (p=1,2),
 \end{align}
where 
$\rk(Q_1)\leq 2m$,  $\rk(Q_2)\leq 2n$ and the matrices $A_p$ are discrete analogs of integration operators:
 \begin{align}\label{T3}&
A_p=\{\cla_{p,i-k}\}_{i,k=1}^n, \quad \cla_{1,r}=\left\{\begin{array}{l}0 \,\, {\mathrm{for}}\,\,  r<0, \\
\I / 2 \,  I_m \,\, {\mathrm{for}}\,\,  r=0, \\ \I \, I_m, \,\, {\mathrm{for}}\,\,  r>0; \end{array} \right.
\\ & \cla_{2,r}=0 \,\, {\mathrm{for}}\,\,  r\not=0, \quad \cla_{2,0}=\{a_{j-\ell}\}_{j,\ell=1}^m;
\quad a_r=\left\{\begin{array}{l}0 \,\, {\mathrm{for}}\,\,  r<0, \\
\I / 2  \,\, {\mathrm{for}}\,\,  r=0, \\ \I , \,\, {\mathrm{for}}\,\,  r>0. \end{array} \right.
\label{T4}
 \end{align}
Here, $A_p^*$
stands for the matrix which is the complex conjugate transpose of $A_p$,
$\rk(Q)$ means the rank of the matrix $Q$, $\I$ is the imaginary unit, and $I_m$ is the $m\times m$ identity matrix. 

Multiplying \eqref{T2} by
$T^{-1}$ from the left and from the right, we obtain the identities $T^{-1}A_p-A_p^*T^{-1}=T^{-1}Q_pT^{-1}$. Both cases
$p=1$ and $p=2$ of these identities provide a tool to recover $T^{-1}$ from $T^{-1}Q_1T^{-1}$ or $T^{-1}Q_2T^{-1}$,
respectively. Such an approach was fruitfully used for Toeplitz  and block Toeplitz matrices. In the case of TBT-matrices,
one has to use the two identities simultaneously, which  is quite nontrivial.

In order to explain our approach in greater detail, introduce polynomial vector function (vector column) $h(y)=h(y_1, y_2)$:
 \begin{align}\label{T7}
h(y)=\{h_r(y)\}_{r=1}^{mn}, \quad & h_r(y):=y_1^{i-1}y_2^{j-1} \quad {\mathrm{for}}\,\, r=m(i-1)+j, 
 \end{align}
$n \geq i \geq 1, \,\, m \geq j \geq 1$.
Assume that $T$ is invertible. Then, the polynomial 
 \begin{align}\label{T8}&
\rho(y,z)=h(\ov{z})^*T^{-1}h(y),
 \end{align}
where $h(\ov{z})=h(\ov{z_1}, \ov{z_2})$ and $\ov{z_p}$ is the complex conjugate of $z_p$,  uniquely determines
$T^{-1}$ and is itself important in astrophysics and signal processing. Further in the paper we study $\rho(y,z)$.
We will need the following notations.

\begin{Nn} \label{Nn1r} The notation $I$ is often used $($for simplicity$)$ instead of $I_r$ for
the identity matrix of the order $r = mn$. By $\1_r$ we denote the vector column of
the order $r$, where all the entries equal $1$, and we often omit the index and write  $\1$ when $r = mn$.
\end{Nn}
 In Section \ref{Str} we show (see \eqref{T13+}) that 
$\rho(y,z)$ is easily expressed via
 \begin{align}\label{T9}&
\om(\la,\mu)=\1^*(A_1^*-\mu_1 I)^{-1}(A_2^*-\mu_2 I)^{-1}T^{-1}(A_2-\la_2 I)^{-1}(A_1-\la_1 I)^{-1}\1.
 \end{align}
Then, we derive an important representation of $\om$ in the main Theorem \ref{TmMom}.
This representation is based on the {\it minimal} information on $T^{-1}$ contained in the matrix
$g_{12}$ (or $g_{21}$, for both matrices see \eqref{T40}) with $4mn$ entries.
In Section \ref{SecT} we complete the study of the structure of $T^{-1}$
(see Theorem \ref{TmOm}). Some auxiliary results on determinants are derived
in Appendix.

The approach works for other important $2$-D structured matrices.

As usual $\BC$ stands for the complex plane and
 $g^{\tau}$ 
stands for the matrix which is the  transpose of $g$.
By diag we denote a diagonal (or block diagonal) matrix. For instance, $\diag\{d_1, d_2\}=\begin{bmatrix}d_1 &0 \\ 0 & d_2 \end{bmatrix}$.
The notation col (column) stands for a column vector or block vector: 
$\col \begin{bmatrix}d_1 & d_2\end{bmatrix}=\begin{bmatrix}d_1 \\ d_2 \end{bmatrix}$. The notation $I_r$ stands for the $r\times r$ identity matrix and the indices of the
matrices $U_r$, $J_{2m}$ and $J_{2n}$ (introduced in \eqref{T42}, \eqref{T43}, \eqref{T48}) also indicate the orders of the corresponding
matrices. In all other cases the indices of the matrices are not related to the order.  Many notations were explained before in the
text of the Introduction (see, e.g., Notation \ref{Nn1r}).


\section{Representation of the $\rho$-polynomial} \label{Str}
\setcounter{equation}{0}
\paragraph{1.} It is immediate that the matrix $A_1$ and the block diagonal matrix $A_2$ (given by \eqref{T3}
and \eqref{T4}) commute.
One easily derives (see e.g. \cite[p. 452]{SaAJFA}) that
 \begin{align}\label{T5}&
(A_1 -\la_1 I)^{-1}\begin{bmatrix}I_m \\ I_m \\ \ldots \\ I_m \end{bmatrix}=\left(\frac{\I}{2}-\la_1\right)^{-1}\begin{bmatrix}I_m \\ \psi(\la_1) I_m\\ \ldots \\
\psi(\la_1)^{n-1} I_m \end{bmatrix}, \quad \displaystyle{\psi(\la_1):=\frac{\la_1 + \frac{\I}{2}}{\la_1 - \frac{\I}{2}}}.
 \end{align}
 As a special case  of \eqref{T5} we have
  \begin{align}\label{T6}&
(\cla_{2,0} -\la_2 I_m)^{-1} \1_m=\left(\frac{\I}{2}-\la_2 \right)^{-1}\begin{bmatrix}1 \\ \psi(\la_2) \\ \ldots \\
\psi(\la_2)^{m-1}  \end{bmatrix}. 
 \end{align}
Since $A_2=\diag\{\cla_{2,0}, \ldots, \cla_{2,0}\}$, relations \eqref{T7}, \eqref{T5} and \eqref{T6} imply that
 \begin{align}\label{T10}&
\left(\frac{\I}{2}-\la_2\right)\left(\frac{\I}{2}-\la_1\right)(A_2-\la_2 I)^{-1}(A_1-\la_1 I)^{-1}\1=h\big(\psi(\la_1), \psi(\la_2)\big).
 \end{align}
Hence, we  also have 
 \begin{align}\label{T11}&
\left(\frac{\I}{2}+\mu_2\right)\left(\frac{\I}{2}+\mu_1\right)\1^*(A_1^*-\mu_1 I)^{-1}(A_2^*-\mu_2 I)^{-1}=h\big(\psi(\ov{\mu_1}), \psi(\ov{\mu_2})\big)^*.
 \end{align}
Further in the text we assume that $T$ is invertible. It is immediate from
 \eqref{T8}, \eqref{T9}, \eqref{T10} and \eqref{T11} that
 \begin{align}\nn
\rho\Big(\psi(\la_1), \psi(\la_2), \ov{\psi(\ov{\mu_1})}, \ov{\psi(\ov{\mu_2})}\Big)=&\left(\frac{\I}{2}+\mu_2\right)\left(\frac{\I}{2}+\mu_1\right)
\\ \label{T12}&
\times \left(\frac{\I}{2}-\la_2\right)\left(\frac{\I}{2}-\la_1\right)\om(\la,\mu).
 \end{align} 
 Recalling that $\psi$ is given in \eqref{T5}, we easily construct the function inverse to $\psi$
 and rewrite the equalities $y_k=\psi(\la_k)$ and $z_k=\ov{\psi\big(\ov{\mu_k}\big)}$ in the form
\begin{align}\label{T13}&
\la_k=\vp(y_k), \quad \mu_k=-\vp(z_k), \quad \vp(\xi):=\frac{\I}{2} \frac{\xi +1}{\xi -1}.
 \end{align}
 Using \eqref{T13}, we rewrite \eqref{T12} in order to express $\rho(y,z)$ via $\om$:
\begin{align}\label{T13+}
\rho\Big(y, z\Big)=&\left(\frac{\I}{2}-\vp(z_2)\right)\left(\frac{\I}{2}-\vp(z_1)\right)
\\ \nn &
\times \left(\frac{\I}{2}-\vp(y_2)\right)\left(\frac{\I}{2}-\vp(y_1)\right)\om\Big(\vp(y_1),\vp(y_2),-\vp(z_1),-\vp(z_2)\Big).
 \end{align} 
 Thus, a representation of $\rho$ will follow from a representation of $\om$.
\paragraph{2.} Next, we show that the matrix  identities \eqref{T2} hold and consider them in detail. 
Indeed, according to \cite[(1.2)]{SaAJFA} we have
 \begin{align}\label{T14}&
A_1T-TA_1^*=\I\big(M_{11}M_{21}+ M_{31}M_{41}\big),
 \end{align}
where 
 \begin{align}\label{T15}&
M_{11}=\col  \begin{bmatrix} \clm_{11}^{(1)} & \clm_{11}^{(2)} & \ldots & \clm_{11}^{(n)}\end{bmatrix}, \quad 
\clm_{11}^{(i)}:=\frac{1}{2}\clt_0 +\sum_{s=1}^{i-1}\clt_{s},
\\ \label{T16}&
M_{21}= \begin{bmatrix} I_m & I_m & \ldots & I_m\end{bmatrix}, \quad M_{31}=M_{21}^*,
\\ \label{T17}&
M_{41}=  \begin{bmatrix} \clm_{41}^{(1)} & \clm_{41}^{(2)} & \ldots & \clm_{41}^{(n)}\end{bmatrix}, \quad 
\clm_{41}^{(k)}:=\frac{1}{2}\clt_0 +\sum_{s=1}^{k-1}\clt_{-s};
\end{align}
We derive the expressions for 
 $\cla_{2,0}\clt_r-\clt_r \cla_{2,0}^*$  as the special  cases of \eqref{T14} (after setting in \eqref{T14}  $m=1$ and substituting then $n=m$).
 In view of these expressions,  for the matrix $A_2=\diag\{\cla_{2,0}, \ldots, \cla_{2,0}\}$ the next identity follows:
  \begin{align}\label{T18}&
A_2T-TA_2^*=\I\big(M_{12}M_{22}+ M_{32}M_{42}\big),
 \end{align}
where $M_{12}$, $M_{32}$ and $M_{22}$, $M_{42}$ are $mn \times n$ and $n \times mn$, respectively, matrices,
 \begin{align}\label{T19}&
M_{12}=\{\clm_{12}^{(i-k)}\}_{i,k=1}^n, \quad M_{22}=\{\clm_{22}^{(i-k)}\}_{i,k=1}^n, 
\\ \label{T20}&
  M_{32}=M_{22}^*, \quad 
M_{42}=\{\clm_{42}^{(i-k)}\}_{i,k=1}^n,
\\ \label{T21}&
 \clm_{12}^{(r)}:=\col \begin{bmatrix} \frac{1}{2}t_r^{(0)} &  \frac{t_r^{(0)}}{2}+t_r^{(1)}\ldots & \frac{t_r^{(0)}}{2}+\sum_{j=2}^{m}t_r^{(j-1)}\end{bmatrix},
\\ \label{T22}&
\clm_{22}^{(r)}=\1_m^* \quad {\mathrm{for}} \quad r=0, \quad \clm_{22}^{(r)}=0 \quad {\mathrm{for}} \quad r\not= 0,
\\ \label{T23}&
\clm_{42}^{(r)}:=\begin{bmatrix} \frac{1}{2}t_r^{(0)} &  \frac{t_r^{(0)}}{2}+t_r^{(-1)}\ldots & \frac{t_r^{(0)}}{2}+\sum_{\ell=2}^{m}t_r^{(1-\ell)}\end{bmatrix}.
\end{align}
Clearly, the matrices $M_{1p}$ and $M_{4p}$ satisfy some operator identities which are similar to
the identities  \eqref{T14} and \eqref{T18}. Indeed, put
 \begin{align}\label{T24}&
\cla_1= \{a_{j-\ell}\}_{j,\ell=1}^n,
\quad a_r=\left\{\begin{array}{l}0 \,\, {\mathrm{for}}\,\,  r<0, \\
\I / 2  \,\, {\mathrm{for}}\,\,  r=0, \\ \I , \,\, {\mathrm{for}}\,\,  r>0; \end{array} \right. \quad \cla_2=\cla_{2,0},
\end{align}
where $\cla_{2,0}$ is given in \eqref{T4}. We note that $\cla_1$ coincides with the special case of $A_1$ where $m=1$,  
and $\cla_2$ coincides with  the special case of $A_2$ where $n=1$. Similarly to \eqref{T14} and \eqref{T18}, one can show
that
 \begin{align}\label{T25}&
\cla_sM_{4p}-M_{4p}A_s^*=\I\clq_s \quad {\mathrm{for}} \,\, s=1, \,\, p=2  \quad  {\mathrm{and}} \quad s=2, \,\, p=1;
\\ \label{T26}&
 \clq_1=K_{11}M_{21}+ \1_n K,  \quad  \clq_2=K_{12}M_{22}+ \1_m K;
 \\ \label{T27}&
 K=\begin{bmatrix}
\frac{1}{2}\clm_{42}^{(0)} &  \, \frac{1}{2}\clm_{42}^{(0)}+\clm_{42}^{(-1)} & \ldots &  \, \frac{1}{2}\clm_{42}^{(0)}+\sum_{k=2}^n\clm_{42}^{(1-k)}
 \end{bmatrix},
  \\ \label{T28}&
 K_{11}=\col \begin{bmatrix}
\frac{1}{2}\clm_{42}^{(0)} &  \frac{1}{2}\clm_{42}^{(0)}+\clm_{42}^{(1)} &\ldots &  \, \frac{1}{2}\clm_{42}^{(0)}+\sum_{i=2}^n\clm_{42}^{(i-1)}
 \end{bmatrix},
 \\ \label{T29}&
 K_{12}=\begin{bmatrix}
\frac{1}{2}\clm_{12}^{(0)} &  \, \frac{1}{2}\clm_{12}^{(0)}+\clm_{12}^{(-1)} &\ldots &  \, \frac{1}{2}\clm_{12}^{(0)}+\sum_{k=2}^n\clm_{12}^{(1-k)}
 \end{bmatrix}.
 \end{align}
 
 \paragraph{3.} The identities \eqref{T14} and \eqref{T18} may be rewritten in the form
  \begin{align}\label{T30}&
A_pT-TA_p^*=\I \Pi_p \wh \Pi_p, \quad \Pi_p:=\begin{bmatrix}M_{1p} & M_{3p}\end{bmatrix}, \quad
\wh \Pi_p:=\begin{bmatrix}M_{2p} \\ M_{4p}\end{bmatrix}.
 \end{align}
Introduce the notations: 
   \begin{align}\label{T31}&
R:=T^{-1}; \quad \G_p:=R\Pi_p, \quad \wh \G_p:=\wh \Pi_p R \quad (p=1,2); 
\\ \label{T31+}&
\G:=\begin{bmatrix}\G_1 & \G_2\end{bmatrix}, \quad 
\wh \G:=\begin{bmatrix}\wh \G_1 \\ \wh \G_2\end{bmatrix}, \quad P_1:\begin{bmatrix}I_{2m} & 0 \\ 0 & 0\end{bmatrix},
\quad P_2:= I_{2(m+n)}-P_1.
 \end{align}
 We start with the following simple representations of $\om$.
 \begin{Pn} \label{simprep} Let $\om(\la, \mu)$ be given by \eqref{T9}. Then we have
  \begin{align}\label{T32}&
\om(\la, \mu)= \I (\la_p-\mu_p)^{-1} u(\mu) P_p \wh u(\la) \quad (p=1,2),
 \end{align}
 where
  \begin{align}\nn
  u(\la)=&\1^*(A_1^*-\mu_1 I)^{-1}(A_2^*-\mu_2 I)^{-1}\G
  \\ & \label{T33}
  -\I \begin{bmatrix}\1_m^*(\cla_2^* -\mu_2 I_m)^{-1} & 0
  & \1_n^*(\cla_1^* -\mu_1 I_n)^{-1} &0
  \end{bmatrix},
 \\ \nn 
\wh u(\la)=& \wh \G (A_2-\la_2 I)^{-1}(A_1-\la_1 I)^{-1}\1
\\ & \label{T34}
+\I \, \col \begin{bmatrix}0 & (\cla_2 -\la_2 I_m)^{-1}\1_m & 0
  & (\cla_1 -\la_1 I_n)^{-1}\1_n 
  \end{bmatrix}.
 \end{align}
 \end{Pn}
\begin{proof}. Relations \eqref{T30} and \eqref{T31} yield
  \begin{align}\label{T35}&
RA_p-A_p^*R=\I \G_p \wh \G_p.
 \end{align}
 Hence, we have 
 \begin{align}\nn&
(\mu_p-\la_p)R= R(A_p-\la_pI)-(A_p^*-\mu_p I)R-\I \G_p \wh \G_p,
 \end{align}
or, equivalently,
 \begin{align}\nn
 (A_p^*-\mu_p I)^{-1}R(A_p-\la_pI)^{-1}= &(\mu_p-\la_p)^{-1}\big((A_p^*-\mu_p I)^{-1}R-R(A_p-\la_pI)^{-1}
\\ \label{T36} & 
-\I (A_p^*-\mu_p I)^{-1}\G_p \wh \G_p(A_p-\la_pI)^{-1}\big).
 \end{align}
 It easily follows from the definitions \eqref{T16}, \eqref{T19} and \eqref{T22} of $M_{2p}$  that
  \begin{align}\label{T37}&
\1^*(A_1^*-\mu_1 I)^{-1}=\1_m^*M_{21}(A_1^*-\mu_1 I)^{-1}=\1_n^*(\cla_1^*- \mu_1 I_n)^{-1}M_{22},
\\ \label{T38}&
\1^*(A_2^*-\mu_2 I)^{-1}=\1_m^*M_{21}(A_2^*-\mu_2 I)^{-1}=\1_m^*(\cla_2^*-\mu_2 I_m)^{-1}M_{21}.
 \end{align}
 Now, set in \eqref{T36} $p=1$, multiply both sides of \eqref{T36} by $\1^*(A_2^*-\mu_2 I)^{-1}$ from
 the left and by $(A_2-\mu_2 I)^{-1}\1$ from the right, and take into account that $A_1$ and $A_2$
 commute and that the equalities \eqref{T9}, \eqref{T38}
 and $M_{31}=M_{21}^*$ hold. Then, we obtain
  \begin{align}\nn
\om(\la,\mu)=& \I(\la_1-\mu_1)^{-1}\Big(\1^*(A_1^*-\mu_1 I)^{-1}(A_2^*-\mu_2 I)^{-1}\G_1 \wh \G_1(A_2-\la_2 I)^{-1}
\\ & \label{T39}
\times (A_1-\la_1 I)^{-1}\1
\\ & \nn
+ \I \1^*(A_1^*-\mu_1 I)^{-1}(A_2^*-\mu_2 I)^{-1}RM_{31}
 (\cla_2-\la_2 I_m)^{-1}\1_m 
 \\ & \nn
 -\I \1_m^*(\cla_2^*-\mu_2 I_m)^{-1}M_{21}R
(A_2-\la_2 I)^{-1}
(A_1-\la_1 I)^{-1}\1\Big).
 \end{align}
Recall the definitions of $\Pi_p$, $\wh \Pi_p$, $\G_p$ and $\wh \G_p$  in \eqref{T30} and \eqref{T31}.
Hence, formula \eqref{T32} follows (for $p=1$) from \eqref{T33}, \eqref{T34} and \eqref{T39}.

Formula \eqref{T32} for $p=2$ follows in the same way (as \eqref{T32} for $p=1$) but this time from \eqref{T36} (with $p=2$)
and from \eqref{T37}.
\end{proof}  
 \paragraph{4.} 
In Proposition \ref{simprep}, we recover $\om$ (and so $\rho$ and $T^{-1}$) either from $\G_1=T^{-1}\Pi_1$ and
$\wh \G_1=\wh \Pi_1T^{-1}$ or from $\G_2=T^{-1}\Pi_2$ and $\wh \G_2=\wh \Pi_2T^{-1}$, which means that only
one of the operator identities (either \eqref{T14} or \eqref{T18}) is used. In this paragraph we will recover $\om$
either  from the  $2m \times 2n$ matrix $\wh \Pi_1 T^{-1} \Pi_2$ or from the $2n \times 2m$ matrix  
$\wh \Pi_2 T^{-1} \Pi_1$, using both operator identities.

First introduce matrices
\begin{align}\label{T40}&
g_{12}=\I \wh \Pi_1 T^{-1} \Pi_2-\I \begin{bmatrix}\1_m \1_n^* & 0 \\ K_{12} & 0\end{bmatrix},
\quad
g_{21}=\I \wh \Pi_2 T^{-1} \Pi_1-\I \begin{bmatrix}\1_n \1_m^* & 0 \\ K_{11} & 0\end{bmatrix},
\end{align}
where $K_{11}$ and $K_{12}$ are given in \eqref{T28} and \eqref{T29}, respectively. The connection
between $g_{12}$ and $g_{21}$ is described by the simple formula
\begin{align}\label{T41}&
g_{21}=-U_{2n} J_{2n}g_{12}^{\tau} J_{2m} U_{2m},
 \end{align}
 where $g_{12}^{\tau} $ is the transpose of $g_{12}$ and
\begin{align}\label{T42}&
U_{2n}=\{\delta_{2n-i-k+1}\}_{i,k=1}^{2n}, 
\quad J_{2n}=\diag\{I_n, -I_n\}, 
\\ \label{T43}&
U_{2m}=\{\delta_{2m-i-k+1}\}_{i,k=1}^{2m}, \quad J_{2m}=\diag\{I_m, -I_m\}.
 \end{align}  
 The validity of \eqref{T41} is shown in the proof of the main Theorem \ref{TmMom} and means
 that given $g_{12}$ we easily construct  $g_{21}$ and vice versa.
 
 Next, we introduce the matrix function
 \begin{align}\label{T44}&
 \sbox0{$\begin{matrix}\cla_2-\la_2 I_m & 0\\ 0 & \cla_2-\la_2 I_m\end{matrix}$}
\sbox1{$\begin{matrix} \cla_1-\la_1 I_n & 0\\ 0 & \cla_1-\la_1 I_n\end{matrix}$}
G(\la):=\left[
\begin{array}{c|c}
\usebox{0}&\makebox[\wd0]{\large $ g_{12}$}\\
\hline
  \vphantom{\usebox{0}}\makebox[\wd0]{\large $ g_{21}$}&\usebox{1}
\end{array}
\right], 
 \end{align} 
 \begin{Rk}\label{RkG}
We introduce a polynomial $\th(\la)$ via the determinant of $G(\la):$
 \begin{align}\label{T45}&
\th(\la)= \sqrt{\det \big(G(\la)\big)},
 \end{align}
 where the branch of the root in \eqref{T45} is chosen so that $\sqrt{\det \big(G(\la)\big)}$
 is a polynomial such that the coefficient corresponding to the term $\la_1^n\la_2^m$
 equals $-1$.
 The existence of the polynomial $\th(\la)$ and the way to construct it explicitly $($when $G(\la)$ of the form \eqref{T44} is given$)$
 are shown in Lemma \ref{LaSqrt}.
 \end{Rk}
 In the following main theorem, we express $\wh u$ and $u$ given by \eqref{T34} and \eqref{T33}, respectively,
 via $G(\la)$ (that is, via $g_{12}$ or $g_{21}$). Then, one can apply Proposition~\ref{simprep} in order to construct
$\om$ and $\rho$. 
\begin{Tm}\label{TmMom} Let a TBT-matrix $T$ be invertible and let the corresponding matrix $g_{12}$
or $g_{21}$ $($of the form \eqref{T40}$)$ be given. 

Then, $\om(\la,\mu)$ of the form \eqref{T9} admits representations \eqref{T32},
where $u$ and $\wh u$ $($which are introduced in \eqref{T33} and \eqref{T34}$)$ may be recovered
from the relations
\begin{align}\label{T46}&
\wh u(\la)=-\I\left(\la_1 -\frac{\I}{2}\right)^{-n}\left(\la_2 -\frac{\I}{2}\right)^{-m}\th(\la)G(\la)^{-1}\col 
\begin{bmatrix}0 & \1_m & 0 & \1_n\end{bmatrix},
\\ \label{T47}&
u(\mu)=\left(\frac{\mu_1-\frac{\I}{2}}{\mu_1+\frac{\I}{2}}\right)^n\left(\frac{\mu_2-\frac{\I}{2}}{\mu_2+\frac{\I}{2}}\right)^m
\wh u(\mu)^{\tau}\diag\{J_{2m}U_{2m}, \, J_{2n} U_{2n}\},
 \end{align}
 $G$ is expressed in \eqref{T44} via $g_{12}$ or, equivalently, via $g_{21}$ $($using \eqref{T41}$)$, and $\th$ is given in \eqref{T45}.
 \end{Tm}
 \begin{proof}. Step 1.  In this step we prove the auxiliary  equalities  \eqref{T47} and \eqref{T41}. For this purpose, recall
 the definitions \eqref{T42} and \eqref{T43} and also set  
\begin{align}\label{T48}&
U=U_{mn}=\{\delta_{mn-i-k+1}\}_{i,k=1}^{mn}, \quad U_m=\{\delta_{m-i-k+1}\}_{i,k=1}^{m}.
 \end{align}   
 It is easy to see that $UTU=T^{\tau}$, and so
\begin{align}\label{T49}&
URU=R^{\tau} \quad (R=T^{-1}).
 \end{align}   
 Similarly to $UTU=T^{\tau}$ we have 
 $$U_m \clt_r^{\, \tau}=\clt_r U_m, \qquad \clt_r^{\, \tau}U_m=U_m\clt_r ,$$  
 which (in view of the definitions \eqref{T15}--\eqref{T17})
 yields
\begin{align}\label{T50}&
UM_{41}^{\tau}+M_{11}U_m=TM_{31}U_m, \quad M_{11}^{\tau} U +U_m M_{41}=U_m M_{21}T.
 \end{align}   
 Using the definitions \eqref{T19}--\eqref{T23}, we derive (in the same way as \eqref{T50}) the
equalities
\begin{align}\label{T51}&
UM_{42}^{\tau}+M_{12}U_n=TM_{32}U_n, \quad M_{12}^{\tau} U +U_n M_{42}=U_n M_{22}T.
 \end{align} 
 Now, let us show that
\begin{align}\label{T52}&
\wh \G_1^{\tau}=U\G_1U_{2m} J_{2m}+\begin{bmatrix}0 & M_{31}\end{bmatrix}, \quad
\wh \G_2^{\tau}=U\G_2U_{2n} J_{2n}+\begin{bmatrix}0 & M_{32}\end{bmatrix}.
 \end{align} 
Indeed,      \eqref{T30}, \eqref{T31}, \eqref{T49} and \eqref{T50} imply that
\begin{align}\nn
\wh \G_1^{\tau}&= URU\begin{bmatrix}M_{21}^{\tau} & M_{41}^{\tau}\end{bmatrix}=
UR\begin{bmatrix}M_{31}U_m &\,\,  -M_{11}U_m+TM_{31}U_m\end{bmatrix}
\\ \nn &
=
UR\begin{bmatrix}M_{11} & M_{31}\end{bmatrix} U_{2m}J_{2m}+\begin{bmatrix}0 & UM_{31}U_m\end{bmatrix},
 \end{align} 
 and the first equality in \eqref{T52} follows. Here we used the immediate equalities 
 \begin{align}\label{T52+}&
M_{21}^{\tau}=M_{31}, \qquad UM_{31}U_m=M_{31}.
 \end{align} 
 Later we will also need an analog of \eqref{T52+} (for the matrices $M_{s2}$), namely, the equalities
  \begin{align}\label{T52!}&
M_{32}^{\tau}=M_{22}, \qquad U_nM_{22}U=M_{22}.
 \end{align} 
Taking into account \eqref{T51}, we (similarly to the proof of the first equality in \eqref{T52}) derive the second equality in \eqref{T52}.
 For
 \begin{align}\label{T53}&
F(\mu_1,\mu_2):=\1^*(A_1^*-\mu_1 I)^{-1}(A_2^*-\mu_2 I)^{-1},
 \end{align}
 we will need the following relation:
\begin{align}\label{T54} &
F(\mu_1,\mu_2)
=
\left(\frac{\mu_1-\frac{\I}{2}}{\mu_1+\frac{\I}{2}}\right)^n\left(\frac{\mu_2-\frac{\I}{2}}{\mu_2+\frac{\I}{2}}\right)^m
\1^*(A_1^*+\mu_1 I)^{-1}(A_2^*+\mu_2 I)^{-1}U.
 \end{align}  
In order to obtain \eqref{T54}, we take into account \eqref{T11} and rewrite \eqref{T53} in the form
\begin{align}\label{T55}
F(\mu_1,\mu_2)&=\{F_r(\mu_1,\mu_2)\}_{r=1}^{mn}
\\ \nn &
=\left(\mu_1+\frac{\I}{2}\right)^{-1}\left(\mu_2+\frac{\I}{2}\right)^{-1}
\left\{\left(\frac{\mu_1-\frac{\I}{2}}{\mu_1+\frac{\I}{2}}\right)^{i-1}\left(\frac{\mu_2-\frac{\I}{2}}{\mu_2+\frac{\I}{2}}\right)^{j-1}\right\},
 \end{align}  
 where $i$ and $j$ for the $r$th entry (of the row vector function) in the braces above are chosen similarly to the way it is done
 in \eqref{T7}. Using the definition \eqref{T53} of $F$, substituting $-\mu_p$ instead of $\mu_p$ $(p=1,2)$ into \eqref{T55} 
 and applying $U$
 from the right,  we derive
  \begin{align}\label{T56}&
\1^*(A_1^*+\mu_1 I)^{-1}(A_2^*+\mu_2 I)^{-1}U
\\ & \nn
=\left(\mu_1-\frac{\I}{2}\right)^{-1}\left(\mu_2-\frac{\I}{2}\right)^{-1}
\left\{\left(\frac{\mu_1+\frac{\I}{2}}{\mu_1-\frac{\I}{2}}\right)^{n-i}\left(\frac{\mu_2+\frac{\I}{2}}{\mu_2-\frac{\I}{2}}\right)^{m-j}\right\}.
 \end{align}
According to \eqref{T56}, the right hand sides of  \eqref{T54} and \eqref{T55} are equal, and so \eqref{T54} follows from
\eqref{T55}.
 
 Partition $u$ and $\wh u$ into the row and column (respectively) vector blocks:
\begin{align}\label{T57}&
 u=\begin{bmatrix}u_1 & u_2\end{bmatrix}, \quad \wh u=\begin{bmatrix}\wh u_1 \\ \wh u_2\end{bmatrix}
 \quad (u_1, \wh u_1 \in \BC^{2m}, \quad u_2, \wh u_2 \in \BC^{2n}).
 \end{align}   
 Taking into account \eqref{T31+}, \eqref{T33}, \eqref{T34}, \eqref{T52} and \eqref{T54}, we see that 
 \begin{align}\label{T58}&
 u_1(\mu)=\left(\frac{\mu_1-\frac{\I}{2}}{\mu_1+\frac{\I}{2}}\right)^n\left(\frac{\mu_2-\frac{\I}{2}}{\mu_2+\frac{\I}{2}}\right)^m
 \wh u_1(\mu)^{\tau}J_{2m}U_{2m}.
  \end{align} 
Here, we also used the relation
\begin{align}\label{T58'}&
\1^*(A_1^*-\mu_1 I)^{-1}(A_2^*-\mu_2 I)^{-1}M_{31}
\\ & \nn
-\I \left(\frac{\mu_1-\frac{\I}{2}}{\mu_1+\frac{\I}{2}}\right)^n
\left(\frac{\mu_2-\frac{\I}{2}}{\mu_2+\frac{\I}{2}}\right)^m \1_m^*(\cla_2^*+\mu_2 I_m)^{-1} U_m
=\I \, \1_m^*(\cla_2^*-\mu_2 I_m)^{-1},
\end{align}   
which easily follows from \eqref{T5} and \eqref{T6}. The next equality, that is,
\begin{align} 
 \label{T59}&
 u_2(\mu)=\left(\frac{\mu_1-\frac{\I}{2}}{\mu_1+\frac{\I}{2}}\right)^n\left(\frac{\mu_2-\frac{\I}{2}}{\mu_2+\frac{\I}{2}}\right)^m
 \wh u_2(\mu)^{\tau}J_{2n}U_{2n}
 \end{align} 
is proved in the same way as \eqref{T58}. Finally, relations \eqref{T57}, \eqref{T58} and \eqref{T59} yield \eqref{T47}.

Relations \eqref{T50} and \eqref{T51} also help  to prove \eqref{T41}.  Indeed, from the definitions \eqref{T40} and \eqref{T30}
we have
\begin{align}\nn
U_{2n} J_{2n}g_{12}^{\tau} J_{2m} U_{2m}=&\I \begin{bmatrix}U_nM_{32}^{\tau} \\ -U_nM_{12}^{\tau} \end{bmatrix}
UT^{-1}U\begin{bmatrix}M_{41}^{\tau}U_m &  -M_{21}^{\tau} U_m\end{bmatrix}
\\ & \label{T60}
-\I
\begin{bmatrix}0&  0 \\ -U_n K_{12}^{\tau} U_m & \1_n \1_m^* \end{bmatrix}.
  \end{align} 
 Using the equalities \eqref{T50} and \eqref{T51} (as well as \eqref{T52+} and \eqref{T52!}), we rewrite \eqref{T60}
 in the form
 \begin{align}& \label{T61}
U_{2n} J_{2n}g_{12}^{\tau} J_{2m} U_{2m}=-\I \begin{bmatrix}M_{22} \\ M_{42}\end{bmatrix}
T^{-1}\begin{bmatrix}M_{11} &  M_{31}\end{bmatrix} +Z,
\end{align} 
where
 \begin{align}\nn
Z=&\I \begin{bmatrix}M_{22} \\ M_{42}\end{bmatrix}
T^{-1}\begin{bmatrix} TM_{31} &  0\end{bmatrix} +\I \begin{bmatrix}0 \\ M_{22}T\end{bmatrix}
T^{-1}\begin{bmatrix} M_{11} &  M_{31}\end{bmatrix}
\\ &  \label{T62}
-\I \begin{bmatrix}0 \\ M_{22}T\end{bmatrix}
T^{-1}\begin{bmatrix} T M_{31} &  0\end{bmatrix}
-\I \begin{bmatrix} 0 & 0 \\ -U_n K_{12}^{\tau} U_m & \1_n\1_m^*\end{bmatrix}
\end{align} 
We partition $Z$ into for blocks $Z=\{Z_{ps}\}_{p,s=1}^2$ and easily derive
 \begin{align}& \label{T63}
Z_{11}=\I M_{22}M_{31}=\I \1_n \1_m^*, \quad Z_{22}=\I(M_{22}M_{31}-\I \1_n \1_m^*)=0,  
\\ & \label{T64}
Z_{12}=0, \quad Z_{21}=\I\big(M_{42}M_{31}+M_{22}(M_{11}-TM_{31})+U_nK_{12}^{\tau}U_m\big).
\end{align} 
In view of \eqref{T52+}, taking the transpose of the second equality in \eqref{T51} we obtain
$M_{11}-TM_{31}=-UM_{41}^{\tau}U_m$. Thus, $Z_{21}$ admits representation
 \begin{align}& \label{T65}
 Z_{21}=\I\big(M_{42}M_{31}-M_{22}UM_{41}^{\tau}U_m+U_nK_{12}^{\tau}U_m\big).
\end{align}
Direct calculations show that
\begin{align}& \label{T66}
 M_{42}M_{31}-M_{22}UM_{41}^{\tau}U_m+U_nK_{12}^{\tau}U_m=K_{11},
\end{align}
where $K_{11}$ is introduced in \eqref{T28}. According to \eqref{T65} and \eqref{T66},
one may rewrite \eqref{T64} as
 \begin{align}& \label{T67}
Z_{12}=0, \quad Z_{21}=\I K_{11}.
\end{align}
Clearly, the definition of $g_{21}$ in \eqref{T40} and relations \eqref{T61}, \eqref{T63} and \eqref{T67}
imply the equality \eqref{T41}.

Step 2. In the Step 2 we prove the basic for our proof equality
\begin{align}\label{T68}&
G(\la)\wh u(\la)=-\I\left(\la_1 -\frac{\I}{2}\right)^{-n}\left(\la_2 -\frac{\I}{2}\right)^{-m}\th(\la)\col 
\begin{bmatrix}0 & \1_m & 0 & \1_n\end{bmatrix},
 \end{align}
 which is equivalent to \eqref{T46}.  For this purpose, we rewrite \eqref{T35} in the form
  \begin{align}& \label{T69}
\G_p\wh \G_p=\I\big((A_p^*-\la_p I)R-R(A_p-\la_p I)\big) \quad (p=1,2).
\end{align}
It is easy to see that 
 \begin{align}
& \label{T74}
 (\cla_1^*-\la_1 I_n)^{-1}M_{22}=M_{22}(A_1^*-\la_1 I)^{-1}, \\
 & \label{T75}
   (\cla_2^*-\la_2 I_m)^{-1}M_{21}=M_{21}(A_2^*-\la_2 I)^{-1},
 \\ & \label{T75+}  
 \1_m^* M_{21}=\1_n^* M_{22}=\1^*.
\end{align} 
Since $M_{3p}=M_{2p}^*$, equalities \eqref{T74}--\eqref{T75+} imply that
  \begin{align}& \label{T70}
M_{31}(\cla_2-\la_2 I_m)^{-1}\1_m=(A_2-\la_2 I)^{-1}\1, \\
& \label{T71}
  M_{32}(\cla_1-\la_1 I_n)^{-1}\1_n=(A_1-\la_1 I)^{-1}\1.
\end{align}
In view of the definitions \eqref{T30}--\eqref{T31+} of the matrices $\wh G$, $\G_p$ and $\Pi_p$,
and in view of the definition  of $\wh u_p(\la)$ (see \eqref{T34} and \eqref{T57}), formulas  \eqref{T69},
\eqref{T70} and \eqref{T71}
yield
\begin{align}\label{T72}&
(A_p^*-\la_p I)^{-1}\G_p \wh u_p(\la)=\I R(A_2-\la_2 I)^{-1}(A_1-\la_1 I)^{-1}\1.
 \end{align}
 Indeed, in the case $p=1$ we have 
 \begin{align}\nn&
(A_1^*-\la_1 I)^{-1}\G_1 \wh u_1(\la)=\I R(A_2-\la_2 I)^{-1}(A_1-\la_1 I)^{-1}\1
\\ & \nn
-\I (A_1^*-\la_1 I)^{-1}R
(A_2-\la_2 I)^{-1}\1+\I(A_1^*-\la_1 I)^{-1}RM_{31}(\cla_2-\la_2 I_m)^{-1}\1_m,
 \end{align}
 and \eqref{T72} follows. In the same way \eqref{T72} is derived for $p=2$.
 
 Taking into account \eqref{T34}, \eqref{T57} and equality $\wh \Pi_p R=\wh \G_p$, 
 we see that \eqref{T72} implies the formula
\begin{align}\label{T73}&
\begin{bmatrix} \mathfrak{A}_2(\la_2) & \I \mathfrak{A}_2(\la_2) \wh \Pi_1 (A_2^*-\la_2 I)^{-1}\G_2
\\
\I \mathfrak{A}_1(\la_1) \wh \Pi_2 (A_1^*-\la_1 I)^{-1}\G_1 &\mathfrak{A}_1(\la_1)
 \end{bmatrix}
\wh u(\la) =\I\begin{bmatrix}0 \\ \1_m \\ 0\\ \1_n
\end{bmatrix},
 \end{align} 
where 
\begin{align}\nn&
\mathfrak{A}_1(\la_1)=\begin{bmatrix}\cla_1^*-\la_1 I_n & 0\\ 0 & \cla_1-\la_1 I_n\end{bmatrix},
\quad
\mathfrak{A}_2(\la_2)=\begin{bmatrix}\cla_2^*-\la_2 I_m & 0\\ 0 & \cla_2-\la_2 I_m\end{bmatrix}.
 \end{align}
Using the last equality in \eqref{T30} and formulas \eqref{T25}, \eqref{T74} and \eqref{T75}, we obtain
the following equalities (for $p=1, \,\, s=2$ and for $p=2, \,\, s=1$): 
\begin{align}\label{T76}&
\wh \Pi_p (A_s^*-\la_s I)^{-1}=\mathfrak{A}_s(\la_s)^{-1}\left(\wh \Pi_p+\I
\begin{bmatrix}0 \\ Q_s\end{bmatrix} (A_s^*-\la_s I)^{-1}\right).
 \end{align}
 Moreover, using again \eqref{T72}  we derive
\begin{align}\label{T77}&
\begin{bmatrix}0 \\ Q_s\end{bmatrix} (A_s^*-\la_s I)^{-1}\G_s\wh u_s(\la)=
\I\begin{bmatrix}0 \\ Q_s\end{bmatrix} R(A_2-\la_2 I)^{-1}(A_1-\la_1 I)^{-1}\1.
 \end{align} 
 Since $Q_s$ is given by \eqref{T26} (where $K$ is a row) and since $M_{2s}R$
 equals $\begin{bmatrix}I_m & 0\end{bmatrix}\wh \G_1$ or  $\begin{bmatrix}I_n & 0\end{bmatrix}\wh \G_2$
 (depending on $s$), one may rewrite the right hand side of \eqref{T77}:
\begin{align}\label{T78}&
\begin{bmatrix}0 \\ Q_1\end{bmatrix} (A_1^*-\la_1 I)^{-1}\G_1\wh u_1(\la)=
\I\begin{bmatrix}0 & 0\\ K_{11} & 0\end{bmatrix} \wh u_1(\la)+\I \wt \theta(\la)\begin{bmatrix}0 \\ \1_n\end{bmatrix} 
\\ \label{T79}&
\begin{bmatrix}0 \\ Q_2\end{bmatrix} (A_2^*-\la_1 I)^{-1}\G_2\wh u_2(\la)=
\I\begin{bmatrix}0 & 0\\ K_{12} & 0\end{bmatrix} \wh u_2(\la)+\I \wt \theta(\la)\begin{bmatrix}0 \\ \1_m\end{bmatrix} ,
\\ \label{T80}&
\wt \theta(\la):=  KR(A_2-\la_2 I)^{-1}(A_1-\la_1 I)^{-1}\1.
 \end{align}  
 In particular, we again used   above the definition \eqref{T34} of $\wh u$. 
 
 Formulas \eqref{T76}, \eqref{T78}
 and \eqref{T79} yield
 \begin{align}\nn
\I \mathfrak{A}_2(\la_2) \wh \Pi_1 (A_2^*-\la_2 I)^{-1}\G_2 \wh u_2(\la)= &\I\left( \wh \Pi_1 \G_2 -
\begin{bmatrix}0 & 0\\ K_{12} & 0\end{bmatrix}\right)
\wh u_2(\la)
\\ \label{T81}&
-\I \wt \theta(\la)\col \begin{bmatrix}0 & \1_m\end{bmatrix},
\\ \nn
\I \mathfrak{A}_1(\la_1) \wh \Pi_2 (A_1^*-\la_1 I)^{-1}\G_1 \wh u_1(\la)= &
\I \left( \wh \Pi_2 \G_1 -
\begin{bmatrix}0 & 0\\ K_{11} & 0\end{bmatrix}\right)
\wh u_1(\la)
\\ \label{T82}& 
- \I \wt \theta(\la)\col \begin{bmatrix}0 & \1_n\end{bmatrix}.
 \end{align}
It follows  from \eqref{T73}, \eqref{T81} and \eqref{T82} that
\begin{align}& \label{T83}
\wt G(\la)\wh u(\la)=\I\big(1+\wt \theta(\la)\big)\col \begin{bmatrix}0 & \1_m &0 & \1_n\end{bmatrix}; \quad 
\wt G(\la)=\{\wt g_{ps}(\la)\}_{p,s=1}^2,
\\ & \label{T84}
 \wt g_{11}(\la):=\mathfrak{A}_2(\la_2), \quad \wt g_{12}(\la)\equiv \wt g_{12}:=\I \wh \Pi_1 \G_2 -
\I \begin{bmatrix}0 & 0\\ K_{12} & 0\end{bmatrix},
\\ & \label{T85}
\wt g_{21}(\la)\equiv \wt g_{21}:=\I \wh \Pi_2 \G_1 -
\I \begin{bmatrix}0 & 0\\ K_{11} & 0\end{bmatrix}, \quad
\wt g_{22}(\la):=\mathfrak{A}_1(\la_1).
\end{align}

According to formulas \eqref{T40}, \eqref{T44}, \eqref{T84} and \eqref{T85}, we have
\begin{align}\label{T86}&
G(\la)=\wt G(\la)+\I \, \col \begin{bmatrix}\1_m & 0 & -\1_n & 0\end{bmatrix}
\begin{bmatrix}\1_m^* & 0 & -\1_n^* & 0\end{bmatrix}.
 \end{align}
 Taking into account the definitions in \eqref{T30}--\eqref{T31+}, the definition of $\wh u$ in \eqref{T34}
 and equalities \eqref{T75+}, we obtain
 \begin{align}\label{T87}&
\begin{bmatrix}\1_m^* & 0 & -\1_n^* & 0\end{bmatrix}\wh \G=0, \quad 
\begin{bmatrix}\1_m^* & 0 & -\1_n^* & 0\end{bmatrix}\wh u(\la)\equiv 0.
 \end{align}
 Relations \eqref{T83}, \eqref{T86} and \eqref{T87} imply that
 \begin{align}& \label{T83'}
G(\la)\wh u(\la)=\I\big(1+\wt \theta(\la)\big)\col \begin{bmatrix}0 & \1_m &0 & \1_n\end{bmatrix}.
\end{align}
Moreover, it is shown in the Appendix (see Lemma \ref{LaTh}) that
 \begin{align}\label{T88}&
\left(\la_1 -\frac{\I}{2}\right)^{-n}\left(\la_2 -\frac{\I}{2}\right)^{-m}\theta(\la)=-\big(1+\wt \theta(\la)\big).
 \end{align}
 Equalities \eqref{T83'} and \eqref{T88} yield \eqref{T68},
 which proves the theorem.
 \end{proof}
\section{The structure of the operator $T^{-1}$} \label{SecT}
\setcounter{equation}{0}
In the previous section we have shown (see Theorem \ref{TmMom}) that given an invertible TBT-matrix
we can recover $\rho$ and $T^{-1}$ from the matrix $g_{12}$ or from $g_{21}$ of the form \eqref{T40}. Here, we complete the description
of the structure of $R=T^{-1}$. Namely, we prove that given an arbitrary $2m \times 2n$ matrix $g_{12}$ or
an arbitrary $2n \times 2m$ matrix $g_{21}$ and using the same formulas as in Section \ref{Str} we recover some
matrix $R$. Moreover, if this $R$ is invertible, then $T=R^{-1}$ is a TBT-matrix.
 
 First, using \eqref{T41} and \eqref{T44} we construct $G(\la)$.  In view of \eqref{T41} and \eqref{T44}, we easily derive
  \begin{align}\label{S1}&
G(\la)=\diag\{-U_{2m}J_{2m}, \, U_{2n} J_{2n}\}G(\la)^{\tau}\diag\{-J_{2m}U_{2m}, \, J_{2n} U_{2n}\}.
 \end{align}
 Next, we obtain $\wh u$ and $u$ via relations \eqref{T46} and \eqref{T47}, where $\th(\la)=\sqrt{\det\big(G(\la)\big)}$.
 The function $\om$ is expressed in \eqref{T32} via $\wh u$ and $u$.
 \begin{La} \label{LaOm} Let $g_{12}$ or $g_{21}$ be given and let the vector functions $\wh u$ and $u$ be
 constructed via  \eqref{T46} and \eqref{T47} $($using also relations \eqref{T41} and \eqref{T44}$)$.
 
 Then $\om(\la,\mu)$, which we obtain from \eqref{T32}, is the same for both cases $p=1$ and $p=2$.
\end{La}
\begin{proof}. The statement of the lemma is equivalent to the
relation
 \begin{align}\label{S8}&
u(\mu)\big((\la_2 -\mu_2)P_1-(\la_1-\mu_1)P_2\big)\wh u(\la)=0
 \end{align}
  In view of \eqref{T46} and \eqref{T47}, the equality \eqref{S8} is equivalent to the equality
\begin{align}\nn &  
\begin{bmatrix}0 & \1_m^{\tau} & 0 & \1_n^{\tau}\end{bmatrix} \left(G(\mu)^{\tau}\right)^{-1}\diag\{J_{2m}U_{2m}, \, J_{2n} U_{2n}\}
\\ \label{S9}&   \times
\big((\la_2 -\mu_2)P_1-(\la_1-\mu_1)P_2\big)G(\la)^{-1}\col \begin{bmatrix}0 & \1_m & 0 & \1_n\end{bmatrix}=0.
 \end{align} 
 Taking into account \eqref{S1}, we rewrite \eqref{S9} as
 \begin{align}\nn&
\begin{bmatrix}0 & \1_m^{\tau} & 0 & \1_n^{\tau}\end{bmatrix} \diag\{-J_{2m}U_{2m}, \, J_{2n} U_{2n}\}G(\mu)^{-1}
\\ \label{S10}&   \times
\big((\la_2 -\mu_2)P_1+(\la_1-\mu_1)P_2\big)G(\la)^{-1}\col \begin{bmatrix}0 & \1_m & 0 & \1_n\end{bmatrix}=0.
 \end{align} 
Since $G(\mu)-G(\la)= (\la_2 -\mu_2)P_1+(\la_1-\mu_1)P_2$, one may rewrite \eqref{S10} in the form
 \begin{align}\nn&
\begin{bmatrix}0 & \1_m^{\tau} & 0 & \1_n^{\tau}\end{bmatrix} \diag\{-J_{2m}U_{2m}, \, J_{2n} U_{2n}\}
\\ \label{S11}&\times
\big(G(\la)^{-1}- G(\mu)^{-1}\big)\col \begin{bmatrix}0 & \1_m & 0 & \1_n\end{bmatrix}=0.
 \end{align} 
Therefore, it remains to prove that
\begin{align} \label{S12}&
\begin{bmatrix}0 & \1_m^{\tau} & 0 & \1_n^{\tau}\end{bmatrix} \diag\{-J_{2m}U_{2m}, \, J_{2n} U_{2n}\}
G(\nu)^{-1}\col \begin{bmatrix}0 & \1_m & 0 & \1_n\end{bmatrix}=0,
 \end{align} 
and the equalities \eqref{S11} and  \eqref{S8} will follow.  Finally, taking the transpose of the left hand side of \eqref{S12}
and using the relations \eqref{S1} and $U_rJ_r=-J_rU_r$ (for the cases $r=2m$ and $r=2n$) we derive
\begin{align} \label{S13} &
\begin{bmatrix}0 & \1_m^{\tau} & 0 & \1_n^{\tau}\end{bmatrix} \diag\{-J_{2m}U_{2m}, \, J_{2n} U_{2n}\}
G(\nu)^{-1}\col \begin{bmatrix}0 & \1_m & 0 & \1_n\end{bmatrix}
\\ \nn &
=
\begin{bmatrix}0 & \1_m^{\tau} & 0 & \1_n^{\tau}\end{bmatrix} 
\big(G(\nu)^{\tau}\big)^{-1}\diag\{-U_{2m}J_{2m}, \, U_{2n} J_{2n}\}\col \begin{bmatrix}0 & \1_m & 0 & \1_n\end{bmatrix}
\\ \nn &
=-\begin{bmatrix}0 & \1_m^{\tau} & 0 & \1_n^{\tau}\end{bmatrix} \diag\{-J_{2m}U_{2m}, \, J_{2n} U_{2n}\}
G(\nu)^{-1}\col \begin{bmatrix}0 & \1_m & 0 & \1_n\end{bmatrix},
 \end{align} 
which yields \eqref{S12}. 
\end{proof} 
\begin{La}\label{LaT} Let the conditions of Lemma \ref{LaOm} hold. Then $\om(\la, \mu)$
determines a unique matrix $R$ such that
 \begin{align}\label{S14}&
\om(\la,\mu)=\1^*(A_1^*-\mu_1 I)^{-1}(A_2^*-\mu_2 I)^{-1}R(A_2-\la_2 I)^{-1}(A_1-\la_1 I)^{-1}\1.
 \end{align}
\end{La} 
\begin{proof}.
According to  Lemma \ref{LaSqrt}, $\,$ the entries of the matrix function $\th(\la)G(\la)^{-1}$,
 where $\th(\la)=\sqrt{\det\big(G(\la)\big)}$, are polynomials. Moreover,   formula \eqref{T44} yields
   \begin{align}\label{S5}&
\lim_{\la_1\to \infty}G(\la)^{-1}=   \begin{bmatrix}\diag\{(\cla_2-\la_2 I_m)^{-1}, \, (\cla_2-\la_2 I_m)^{-1}\} & 0\\ 0 &0  \end{bmatrix},
\\
\label{S6}&
\lim_{\la_2\to \infty}G(\la)^{-1}=   \begin{bmatrix}0 & 0 \\ 0 &\diag\{(\cla_1-\la_1 I_n)^{-1}, \, (\cla_1-\la_1 I_n)^{-1}\} \end{bmatrix},
 \end{align}
as well as the asymptotic relation
  \begin{align}
\label{S7}&
\th(\la)=-q(\la)\big(1+o(1)\big), \quad q(\la):=\left(\la_1 -\frac{\I}{2}\right)^{n}\left(\la_2 -\frac{\I}{2}\right)^{m},
 \end{align}
when either $\la_1 \to \infty$ or $\la_2 \to \infty$.  Thus,
 for the vector polynomial
   \begin{align}\label{S2}&
v(\la)= \th(\la)G(\la)^{-1} \begin{bmatrix}0 \\ \1_m \\ 0
  \\ \1_n 
  \end{bmatrix}
  +q(\la)  \begin{bmatrix}0 \\ (\cla_2 -\la_2 I_m)^{-1}\1_m \\ 0
  \\ (\cla_1 -\la_1 I_n)^{-1}\1_n 
  \end{bmatrix}
 \end{align}
we have
  \begin{align}\label{S3}&
\lim_{\la_1 \to \infty}\left(\la_1 -\frac{\I}{2}\right)^{-n}v(\la)=0, \quad
\lim_{\la_2 \to \infty}\left(\la_2 -\frac{\I}{2}\right)^{-m}v(\la)=0.
 \end{align}
In other words, the degrees of $\la_1$ in the entries of $v$ are less than $n$ and the degrees of $\la_2$ in the entries of $v$ are less than $m$. 

Similar to the partitioning \eqref{T57} for $\wh u$, partition $v$ into the blocks $v=\col \begin{bmatrix}v_1 & v_2   \end{bmatrix}$.
Using \eqref{T32}, \eqref{T45}, \eqref{T46} and \eqref{S2}, we obtain
  \begin{align}
\label{S15} 
\Om(\la,\mu):&=q(\la)\left(\mu_1 +\frac{\I}{2}\right)^{n}\left(\mu_2 +\frac{\I}{2}\right)^{m}\om(\la,\mu)
\\ \nn &
=-\frac{\I}{\la_1-\mu_1}V_1(\mu)^{\tau}J_{2m}U_{2m}V_1(\la)
=-\frac{\I}{\la_2-\mu_2}V_2(\mu)^{\tau}J_{2n}U_{2n}V_2(\la),
 \end{align}
where $V_1(\la)=v_1(\la)-q(\la) \begin{bmatrix}0 \\ (\cla_2 -\la_2 I_m)^{-1}\1_m   \end{bmatrix}$ and \\
$ V_2(\la)=v_2(\la)-q(\la) \begin{bmatrix}0 \\ (\cla_1 -\la_1 I_n)^{-1}\1_n   \end{bmatrix}.$ By virtue of \eqref{S15}, we have
$$(\la_2-\mu_2)V_1(\mu)^{\tau}J_{2m}U_{2m}V_1(\la)=(\la_1-\mu_1)V_2(\mu)^{\tau}J_{2n}U_{2n}V_2(\la),$$
where $\la_p-\mu_p$ are irreducible polynomials. Hence, $V_1(\mu)^{\tau}J_{2m}U_{2m}V_1(\la)$ may be factored
into the product of polynomials with $\la_1-\mu_1$ as one of the factors (and a similar statement is valid for $V_2(\mu)^{\tau}J_{2n}U_{2n}V_2(\la)$
and $\la_2-\mu_2$).  Thus, we see that $\Om(\la,\mu)$ introduced in \eqref{S15} is a polynomial.

Then, the degrees of $\la_1$ and $\mu_1$ in  the terms of $\Om$ are less than $n$ and the degrees of $\la_2$ and $\mu_2$ are less than $m$ (since the degrees of $\la_1$ in the entries of $v(\la)$ are less than $n$ and the degrees of $\la_2$ in the entries of $v(\la)$ are less than $m$). In order to derive that, we also take the factors $(\la_p - \mu_p)^{-1}$ in \eqref{S15} into account.

It is easy to see that the polynomials 
$$p_s(\nu)=\left(\nu-\frac{\I}{2}\right)^{s-1}\left(\nu +\frac{\I}{2}\right)^{r-s} \quad(1\leq s\leq r)$$ 
are linearly independent polynomials (in one variable $\nu$). 
Therefore, the mentioned above bounds on the degrees of the variables in $\Om(\la,\mu)$ show that
$\Om(\la,\mu)$ may be presented as a linear combination of  the products 
 \begin{align}
\label{S16} 
p_{s_i}(\la_1)p_{s_k}(\mu_1)
p_{s_j}(\la_2)p_{s_{\ell}}(\mu_2) \quad (1 \leq s_i, \,s_k \leq n, \quad 1 \leq s_j, \,s_{\ell} \leq m).
 \end{align}
On the other hand, relations \eqref{T7}, \eqref{T10} and \eqref{T11} (where $\psi$ is introduced in \eqref{T5})
show that the multiplication
of the right hand side of \eqref{S14} by
$q(\la)\left(\mu_1 +\frac{\I}{2}\right)^{n}\left(\mu_2 +\frac{\I}{2}\right)^{m}$ brings us a  polynomial, which
is a linear combination of  the  same products \eqref{S16} with the entries of $R$ as the coefficients.
It is immediate that there is a unique matrix $R$ such that \eqref{S14} holds.
\end{proof}
\begin{Tm} \label{TmOm} Let a $2m \times 2n$ matrix $g_{12}$ or
a  $2n \times 2m$ matrix $g_{21}$ be given and let the vector functions $\wh u$ and $u$ be
 constructed via  \eqref{T46} and \eqref{T47} $($using also relations \eqref{T41} and \eqref{T44}$)$.
 The matrix  $\om(\la,\mu)$ given $($in terms of $u$ and $\wh u)$ by \eqref{T32} determines a unique matrix $R$ such that \eqref{S14} holds.
 Assume that $\det R\not=0$.

Then, the matrix $T=R^{-1}$ is a  TBT-matrix.
\end{Tm}
\begin{proof}.  Step 1. The unique recovery of $R$ from $g_{12}$ or, equivalently, from $g_{21}$ is described in Lemma \ref{LaT}. Now,
consider again $v(\la)$ given by \eqref{S2}.
Since the degrees of $\la_1$ in the entries of $v$ are less than $n$ and the degrees of $\la_2$ in the entries of $v$ are less than $m$,
one can introduce  $\wh \G$ by the equality
 \begin{align}
\label{S17} 
q(\la)\wh \G(A_2-\la_2 I)^{-1}(A_1-\la_1 I)^{-1}\1=-\I v(\la).
 \end{align}
 We partition $\wh \G$ into the blocks $\wh \G_1$ and $\wh \G_2$ (as in \eqref{T31+}), and express matrices
 $\G_p$ ($p=1,2)$ via $\wh \G_p$ using equalities
\begin{equation}\label{S18}
\G_1=U\left(\wh \G_1^{\tau}-\begin{bmatrix}0 & M_{31}\end{bmatrix}\right)J_{2m}U_{2m}, 
\quad \G_2=U\left(\wh \G_2^{\tau}-\begin{bmatrix}0 & M_{32}\end{bmatrix}\right)J_{2n}U_{2n},
 \end{equation} 
which are equivalent to \eqref{T52}. 

Step 2. The basic step in the theorem's proof is the proof of the matrix identities
 \begin{align}&
\label{S19} 
RA_p-A_p^*R=\I \G_p \wh \G_p \quad (p=1,2).
 \end{align}
Clearly,  identities \eqref{S19} are equivalent to the identities 
$$(\mu_p-\la_p)R=R(A_p-\la_pI)-(A_p^*-\mu_p I)R-\I \G_p \wh \G_p.$$
Hence, in view of \eqref{T7}, \eqref{T8} and \eqref{T11},  the identities \eqref{S19} are equivalent to the equalities
 \begin{align}\nn&
\1^*(A_1^*-\mu_1 I)^{-1}(A_2^*-\mu_2 I)^{-1}R(A_2-\la_2 I)^{-1}(A_1-\la_1 I)^{-1}\1
\\ & \nn
=\1^*\big((A_1^*-\mu_1 I)^{-1}(A_2^*-\mu_2 I)^{-1}R(A_s-\la_s I)^{-1}
\\ &  \label{S20} \quad
-(A_s^*-\mu_s I)^{-1}R(A_2-\la_2 I)^{-1}(A_1-\la_1 I)^{-1}
\\ & \nn \quad
-\I
(A_1^*-\mu_1 I)^{-1}(A_2^*-\mu_2 I)^{-1}\G_p\wh \G_p(A_2-\la_2 I)^{-1}(A_1-\la_1 I)^{-1}\big)\1 \big/ (\mu_p-\la_p),
\end{align}
where either $p=1, \, s=2$ or $p=2, \, s=1$. 

Recall that $R$ is determined by $\om(\la,\mu)$ via  \eqref{S14}, where $\om(\la,\mu)$
is given by \eqref{T32}. Setting in \eqref{T32} $p=1$, we prove \eqref{S20} for $p=1$. 
For this purpose, we simplify the terms on the right hand side of  \eqref{S20}.
First, note that
formulas \eqref{T32} and \eqref{S14} imply the equalities
 \begin{align}\nn
\1^*(A_1^*-\mu_1 I)^{-1}(A_2^*-\mu_2 I)^{-1}R(A_2-\la_2 I)^{-1}\1&=\lim_{\la_1\to \infty}\big(-\la_1\om(\la,\mu)\big)
\\ \label{S21}  &
=
-\I u_1(\mu)\lim_{\la_1 \to \infty} \wh u_1(\la).
 \end{align}
According to \eqref{T46} and \eqref{S2} we have
 \begin{align}&
\label{S22} 
\wh u_1(\la)=-\I \left(q(\la)^{-1}v_1(\la)-\begin{bmatrix}0 \\ (\cla_2-\la_2 I_m)^{-1}\1_m \end{bmatrix}\right).
 \end{align} 
Hence, taking into account the first equality in \eqref{S2} and the second equality in \eqref{S7}, we derive
 \begin{align}&
\label{S23} 
\lim_{\la_1 \to \infty} \wh u_1(\la)=\I \begin{bmatrix}0 \\ (\cla_2-\la_2 I_m)^{-1}\1_m \end{bmatrix}.
 \end{align} 
 In view of  \eqref{S23}, we rewrite \eqref{S21} in the form
 \begin{align}\nn &
\1^*(A_1^*-\mu_1 I)^{-1}(A_2^*-\mu_2 I)^{-1}R(A_2-\la_2 I)^{-1}\1
\\ \label{S24}  &
=
 u_1(\mu)\col \begin{bmatrix}0 & (\cla_2-\la_2 I_m)^{-1}\1_m \end{bmatrix}.
 \end{align}
Moreover, relations \eqref{T47} and \eqref{S23} imply that
 \begin{align}\nn
\lim_{\mu_1 \to \infty}  u_1(\mu)&=-\I \left(\mu_2-\frac{\I}{2}\right)^m \left(\mu_2+\frac{\I}{2}\right)^{-m}
\begin{bmatrix}0 & \1_m^{\tau}(\cla_2^{\tau}-\mu_2 I_m)^{-1} \end{bmatrix}U_{2m}
\\ \label{S25}  &
=\I \left(\mu_2-\frac{\I}{2}\right)^m \left(\mu_2+\frac{\I}{2}\right)^{-m}
\begin{bmatrix} \1_m^{*}(\cla_2^{*}+\mu_2 I_m)^{-1} U_m & 0\end{bmatrix}.
 \end{align}
Using \eqref{T6} in order to calculate the right hand side of \eqref{S25}, after easy transformations we have
 \begin{align}\label{S25+}  &
\lim_{\mu_1 \to \infty}  u_1(\mu)=-\I 
\begin{bmatrix} \1_m^{*}(\cla_2^{*}-\mu_2 I_m)^{-1} & 0\end{bmatrix}.
 \end{align}
By virtue of \eqref{T32}, \eqref{S14} and \eqref{S25+}, we obtain a result
 \begin{align} \nn &
 \1^*(A_2^*-\mu_2 I)^{-1}R(A_2-\la_2 I)^{-1}(A_1-\la_1 I)^{-1}\1
\\ \label{S26}  & 
 =\lim_{\mu_1\to \infty}\big(-\mu_1\om(\la,\mu)\big)
=\begin{bmatrix} \1_m^{*}(\cla_2^{*}-\mu_2 I_m)^{-1} & 0\end{bmatrix}\wh u_1(\la).
 \end{align}
 It remains to simplify the expression
\begin{align}\label{S27}  &
\1^*(A_1^*-\mu_1 I)^{-1}(A_2^*-\mu_2 I)^{-1}\G_1\wh \G_1(A_2-\la_2 I)^{-1}(A_1-\la_1 I)^{-1}\1.
 \end{align}
 According to \eqref{S17} and \eqref{S22} we derive
\begin{align}\label{S28}  
\wh \G_1(A_2-\la_2 I)^{-1}(A_1-\la_1 I)^{-1}\1&=-\I q(\la)^{-1}v(\la)
\\ \nn
&=\wh u_1(\la)-\I \,
\col \begin{bmatrix} 0 & (\cla_2-\la_2 I_m)^{-1} \1_m\end{bmatrix}.
 \end{align} 
From \eqref{T53}, \eqref{T54} and \eqref{S18} we also see  that
\begin{align}\nn  &
\1^*(A_1^*-\mu_1 I)^{-1}(A_2^*-\mu_2 I)^{-1}\G_1
\\ \nn &
=
\left(\frac{\mu_1-\frac{\I}{2}}{\mu_1+\frac{\I}{2}}\right)^n\left(\frac{\mu_2-\frac{\I}{2}}{\mu_2+\frac{\I}{2}}\right)^m
\1^*(A_1^{\tau}-\mu_1 I)^{-1}(A_2^{\tau}-\mu_2 I)^{-1}
\\ \label{S29} & \quad \times
\left(\wh \G_1^{\tau}-\begin{bmatrix}0 & M_{31}\end{bmatrix}\right)J_{2m}U_{2m}.
 \end{align}
Furthermore, formula \eqref{T58'} implies the equality
 \begin{align}\label{S30}  &
\1^*(A_1^{\tau}-\mu_1 I)^{-1}(A_2^{\tau}-\mu_2 I)^{-1}M_{31}
\\ & \nn
=\I\1_m^*(\cla_2^*+\mu_2 I_m)^{-1}
+\I \left(\frac{\mu_1+\frac{\I}{2}}{\mu_1-\frac{\I}{2}}\right)^n\left(\frac{\mu_2+\frac{\I}{2}}{\mu_2-\frac{\I}{2}}\right)^m
\1_m^*(\cla_2^*-\mu_2 I_m)^{-1}U_m.
 \end{align}
 Substituting \eqref{S30} into \eqref{S29} and taking into account \eqref{S28}, we derive
\begin{align}\label{S31}  &
\1^*(A_1^*-\mu_1 I)^{-1}(A_2^*-\mu_2 I)^{-1}\G_1
\\ \nn &
=
\left(\frac{\mu_1-\frac{\I}{2}}{\mu_1+\frac{\I}{2}}\right)^n\left(\frac{\mu_2-\frac{\I}{2}}{\mu_2+\frac{\I}{2}}\right)^m
\wh u_1(\mu)^{\tau}J_{2m}U_{2m}+\I \begin{bmatrix}\1_m^*(\cla_2^*-\mu_2 I_m)^{-1}& 0 \end{bmatrix}.
 \end{align}
Using \eqref{T47}, we rewrite \eqref{S31}:
 \begin{align}\label{S32}  &
\1^*(A_1^*-\mu_1 I)^{-1}(A_2^*-\mu_2 I)^{-1}\G_1=
u_1(\mu)+\I \begin{bmatrix}\1_m^*(\cla_2^*-\mu_2 I_m)^{-1}& 0 \end{bmatrix}.
 \end{align}
 Finally, denoting the right hand side of \eqref{S20} (where $p=1$ and $s=2$) by
 $Y(\la,\mu)$ and taking into account \eqref{S24}, \eqref{S26}, \eqref{S28} and \eqref{S32} we obtain
   \begin{align}\nn
Y(\la,\mu)=&\Big( u_1(\mu)\col \begin{bmatrix}0 & (\cla_2-\la_2 I_m)^{-1}\1_m \end{bmatrix}
-\begin{bmatrix} \1_m^{*}(\cla_2^{*}-\mu_2 I_m)^{-1} & 0\end{bmatrix}\wh u_1(\la)
\\ \nn &
-\I \big(
u_1(\mu)+\I \begin{bmatrix}\1_m^*(\cla_2^*-\mu_2 I_m)^{-1}& 0 \end{bmatrix}
\big)
\\ &\label{S33}  \times
\big(\wh u_1(\la)-\I \,
\col \begin{bmatrix} 0 & (\cla_2-\la_2 I_m)^{-1} \1_m\end{bmatrix}\big)
\Big)\Big/(\mu_1-\la_1).
 \end{align}
Collecting similar terms in \eqref{S33} and using \eqref{T32}, we derive
   \begin{align}\label{S34} &
Y(\la,\mu)=\I (\la_1-\mu_1)^{-1}u_1(\mu)
\wh u_1(\la)=\om(\la,\mu).
 \end{align}
That is, \eqref{S20} (and so \eqref{S19}) with $p=1$ is equivalent to \eqref{S14},
which proves \eqref{S19} for $p=1$. In the same way, \eqref{S19} is proved for $p=2$.

Step 3. Next, we prove that
   \begin{align}\label{S35} &
\G_1 \begin{bmatrix} 0 \\ I_m\end{bmatrix}=RM_{31}, \quad  \begin{bmatrix} I_m & 0\end{bmatrix}\wh \G_1=M_{21}R,
 \end{align}
 where $M_{21}$ and $M_{31}$ have the special forms \eqref{T16} and do not depend on $R$.
Indeed, the identity \eqref{S19} for $p=1$ implies that 
$$(A_1^*-\mu_1 I)^{-1}R-R(A-\mu_1I)^{-1}\equiv \I(A_1^*-\mu_1I)^{-1}\G_1\wh \G_1(A-\mu_1I)^{-1},$$
and, in particular, we have
\begin{align}\nn &
M_{21}(A_1^*-\mu_1I)^{-1}RM_{31}-M_{21}R(A-\mu_1I)^{-1}M_{31}
\\ \label{S36} &
\equiv \I M_{21}(A_1^*-\mu_1I)^{-1}\G_1\wh \G_1(A-\mu_1I)^{-1}M_{31}.
 \end{align}
Furthermore, equality \eqref{S8} yields
 \begin{align}\label{S37} &
u_1(\mu)\wh u_1(\la)\equiv 0 \quad {\mathrm{for}} \quad \la_1=\mu_1.
 \end{align}
On the other hand, using \eqref{T38} in order to rewrite the left hand sides of \eqref{S28} and \eqref{S32}, 
we obtain the following expressions
for $\wh u_1$ and $u_1$, respectively:
 \begin{align}\label{S38} &
\wh u_1(\la)=\left(\wh \G_1(A_1-\la_1 I)^{-1}M_{31}+\I \begin{bmatrix} 0 \\ I_m\end{bmatrix}\right)
(\cla_2-\la_2 I_m)^{-1}\1_m,
\\ \label{S39} &
u_1(\mu)=\1_m^*(\cla_2^*-\mu_2 I_m)^{-1}
\left(M_{21}(A_1^*-\mu_1 I)^{-1}\G_1-\I \begin{bmatrix}  I_m & 0\end{bmatrix}\right).
 \end{align}
Substituting \eqref{S38} and \eqref{S39} into \eqref{S37}, we see that
 \begin{align}\label{S40} &
\left(M_{21}(A_1^*-\mu_1 I)^{-1}\G_1-\I \begin{bmatrix}  I_m & 0\end{bmatrix}\right)
\left(\wh \G_1(A_1-\mu_1 I)^{-1}M_{31}+\I \begin{bmatrix} 0 \\ I_m\end{bmatrix}\right)\equiv 0.
 \end{align}
It is immediate from \eqref{S36} and \eqref{S40} that
\begin{align} \nn &
M_{21}(A_1^*-\mu_1I)^{-1}\left(RM_{31}-\G_1 \begin{bmatrix} 0 \\ I_m\end{bmatrix}\right)
\\ \label{S41} &
\equiv \left(M_{21}R- \begin{bmatrix}  I_m  & 0 \end{bmatrix}\wh \G_1\right)(A-\mu_1I)^{-1}M_{31}.
 \end{align}
Taking into account the expressions for  $\, M_{21}(A_1^*-\mu_1I)^{-1}\,$ and  for \\ $(A-\mu_1I)^{-1}M_{31}$
(see \eqref{T5} and \eqref{T16}), we conclude from \eqref{S41} that \eqref{S35} holds.

In a similar way, one proves that
   \begin{align}\label{S42} &
\G_2 \begin{bmatrix} 0 \\ I_n\end{bmatrix}=RM_{32}, \quad  \begin{bmatrix} I_n & 0\end{bmatrix}\wh \G_2=M_{22}R,
 \end{align}
where $M_{22}$ and $M_{32}$ are introduced in \eqref{T19}, \eqref{T20}, \eqref{T22}.

Since we assumed that $R$ is invertible, relations \eqref{S19}, \eqref{S35} and \eqref{S42} yield
the matrix identities
 \begin{align}\label{S43}&
A_pT-TA_p^*=\I T\G_p \wh \G_p T=\I\big(M_{1p}M_{2p}+ M_{3p}M_{4p}\big) \quad (T:=R^{-1}),
 \end{align}
where $M_{2p}$ and $M_{3p}$ are given by \eqref{T16}, \eqref{T19}, \eqref{T20}, \eqref{T22}, and
\begin{align}\nn &
M_{11}=T\G_1 \begin{bmatrix} I_m \\ 0\end{bmatrix},  \quad  M_{12}=T \G_2 \begin{bmatrix} I_n \\ 0 \end{bmatrix},
\\ \nn &
M_{41}=\begin{bmatrix} 0 & I_m\end{bmatrix}\wh \G_1 T, \quad M_{42}=\begin{bmatrix} 0 & I_n\end{bmatrix}\wh \G_2 T.
 \end{align}
In the case $p=1$ (and with $M_{21}$ and $M_{31}$ given by \eqref{T16}), the identity \eqref{S43}
means that $T$ is a block Toeplitz matrix (i.e.,  $T=\{\clt_{i-k}\}_{i,k=1}^n$). Indeed, 
 when the right hand side of \eqref{S43} is fixed, the matrix $T$ satisfying \eqref{S43} is unique,
and there is always a block Toeplitz matrix $T$ determined by  \eqref{S43} with $p=1$. One can construct
this $T$  by rewriting  \eqref{S43} ($p=1$) in the form
$$A_1T-TA_1^*=\I\big(\wt M_{11}M_{21}+ \wt M_{31}M_{41}\big); \,\, \wt M_{11}=M_{11}+\b M_{31},
\,\, \wt M_{41}=M_{41}-\b M_{21},$$
where $\b$ is chosen so that the first blocks of $\wt M_{11}$ and $\wt M_{41}$ coincide.
After that one writes down the blocks of $\wt M_{11}$ and $\wt M_{41}$ in the form of the last equalities
in \eqref{T15} and \eqref{T17}, respectively. In this way, one determines the blocks  $\clt_{i-k}$ of
the block Toeplitz matrix $T$  satisfying \eqref{S43} with $p=1$.
 
Similarly, the identity \eqref{S43}, where $p=2$, means that the blocks $\clt_{ik}$ of $T$
are Toeplitz. Thus, $T$ is TBT-matrix.
\end{proof}

\noindent{\bf Acknowledgments.} 
 {The research   was supported by the
Austrian Science Fund (FWF) under Grant  No. P29177.}
\appendix
\section{Auxiliary results}\label{Pr}
\setcounter{equation}{0}
Here, we study the expression $\sqrt{\det G(\la)}$.
The fact that the  ring $\BC[x_1,...x_n]$ of polynomials (in several variables)  with complex
coefficients is factorial, that is, that $\BC[x_1,...x_n]$ is a unique factorisation domain, is well known (see, e.g., \cite[p. 72]{Sha}),
and it is essential in our
further considerations.
We will also need  some interrelations between the minors of the given matrices and their cofactor matrices (see, e.g., \cite{GaKr}).

First note that in view of \eqref{T41} and \eqref{T44}  we obtain \eqref{S1} (both in Sections \ref{Str} and  \ref{SecT}).
Thus, we have a skew-symmetric matrix function
 \begin{align}&
\label{A1} 
\breve G(\la):=\diag \{J_{2m}U_{2m}, -J_{2n}U_{2n}\}G(\la)=-\breve G(\la)^{\tau}.
 \end{align}
Clearly, \eqref{A1} yields
 \begin{align}
& \label{A2} 
\det\breve G(\la)=\det G(\la).
 \end{align}
\begin{Nn} \label{NnD} The notation $D^{i_1 i_2 \ldots i_s}_{k_1 k_2 \ldots k_s}(\la)$ stands for the determinant of the matrix,
cut down from  $\breve G(\la)$ by removing its rows with the numbers $i_1, i_2 ,\ldots, i_s$  and its columns
with the numbers $k_1, k_2 ,\ldots, k_s$ $(i_j\not=i_{\ell}$ and $k_j\not=k_{\ell}$ for $j\not=\ell)$.

In our notations,  $D(\la)=\det \breve G(\la)$ $($if $s=0)$, and we set $D^{1 \, 2  \, 2m}_{1 \, 2  \, 2m}(\la)=1$.
\end{Nn}
According to the well-known property of the minors  (see, e.g., formula (4) in $\S$2 of Ch. 1 \cite{GaKr} and set there $p=2$),
which we apply to the matrix 
cut down from  $\breve G(\la)$ by removing the rows and columns with the same numbers 
$i_1, i_2, \ldots, i_{r-1}$ ($1 \leq r \leq 2m+2n-1$), we have
\begin{align}&
\nn
D^{i_1  \ldots i_{r-1} i_r}_{i_1  \ldots i_{r-1} i_r}(\la)D^{i_1  \ldots i_{r-1} i_{r+1}}_{i_1  \ldots i_{r-1} i_{r+1}}(\la)
-D^{i_1  \ldots i_{r-1} i_r}_{i_1  \ldots i_{r-1} i_{r+1}}(\la)D^{i_1  \ldots i_{r-1} i_{r+1}}_{i_1  \ldots i_{r-1} i_r}(\la)
\\ \label{A3}  &=
D^{i_1  \ldots i_{r-1} }_{i_1  \ldots i_{r-1} }(\la)D^{i_1  \ldots i_{r-1} i_r i_{r+1}}_{i_1  \ldots i_{r-1} i_r i_{r+1}}(\la).
 \end{align}
 We note that either the numbers of the rows and columns, which are cut down, or the numbers of the
 rows and columns, which are preserved, are given for minors in the literature (and in this respect our notations differ
 from the notations in \cite{GaKr}).
 
The skew-symmetric structure of $\breve G$ (see \eqref{A1}) implies for the odd values of $r$ the equalities
\begin{align}&
\label{A4}
D^{i_1  \ldots i_{r-1} i_r}_{i_1  \ldots i_{r-1} i_r}(\la)=D^{i_1  \ldots i_{r-1} i_{r+1}}_{i_1  \ldots i_{r-1} i_{r+1}}(\la)=0,
\quad
D^{i_1  \ldots i_{r-1} i_r}_{i_1  \ldots i_{r-1} i_{r+1}}(\la)=-D^{i_1  \ldots i_{r-1} i_{r+1}}_{i_1  \ldots i_{r-1} i_r}(\la).
 \end{align}
From  \eqref{A3} and \eqref{A4} we derive (for the odd values of $r$)
an important relation:
\begin{align}&
 \label{A5} 
D^{i_1  \ldots i_{r-1} }_{i_1  \ldots i_{r-1} }(\la)=\big(D^{i_1  \ldots i_{r-1} i_r}_{i_1  \ldots i_{r-1} i_{r+1}}(\la)\big)^2
\big(D^{i_1  \ldots i_{r-1} i_r i_{r+1}}_{i_1  \ldots i_{r-1} i_r i_{r+1}}(\la)\big)^{-1}.
 \end{align}
\begin{La}\label{LaSqrt} Let $G(\la)$ of the form \eqref{T44} be given and let \eqref{T41} hold.
Then, $\sqrt{\det G(\la)}$ is a polynomial, and it is presented by the formula
\begin{align}&
 \label{A6} 
\sqrt{\det G(\la)}=\a\frac{D^{i_1}_{i_2}(\la)  D^{i_1  \ldots i_{4} i_5}_{i_1  \ldots i_{4} i_{6}}(\la)\ldots 
D^{i_1  \ldots i_{r-1} i_r}_{i_1  \ldots i_{r-1} i_{r+1}}(\la)}{D^{i_1 i_2 i_3}_{i_1 i_2 i_4}(\la)  \ldots D^{i_1  \ldots i_{s-1} i_s}_{i_1  \ldots i_{s-1} i_{s+1}}(\la)},
 \end{align}
where $\a= \pm 1$ and we fix $\a$ $($in accordance with Remark \ref{RkG}$)$
so that the coefficient corresponding to $\la_1^n\la_2^m$ in $\sqrt{\det G(\la)}$ equals $-1$. $($Here, $r+3=s+1=2m+2n$
if  $2m+2n$ is divisible by $4$, and $r+1=s+3=2m+2n$ if $2m+2n$ is not.$)$

Moreover,  the entries of $\sqrt{\det G(\la)}\,G(\la)^{-1}$ are polynomials as well.
\end{La}
\begin{proof}. Setting $r=1$ in \eqref{A5} and taking into account Notation \ref{NnD} and equality \eqref{A2}, we obtain
\begin{align}&
\label{A7}
\det G(\la)=\det \breve G(\la)=\big(D^{i_1}_{i_2}(\la)\big)^2 \big/ D^{i_1   i_{2}}_{i_1   i_{2}}(\la). 
\end{align}
In the same way, we set $r=3$ in \eqref{A5} and express $D^{i_1   i_{2}}_{i_1   i_{2}}$ via determinants of  smaller
matrices. Next, we set $r=5, \ldots$, and after a final number of steps we have
\begin{align}&
 \label{A8} 
\det G(\la)=\left(\frac{D^{i_1}_{i_2}(\la)  D^{i_1  \ldots i_{4} i_5}_{i_1  \ldots i_{4} i_{6}}(\la)\ldots 
D^{i_1  \ldots i_{r-1} i_r}_{i_1  \ldots i_{r-1} i_{r+1}}(\la)}{D^{i_1 i_2 i_3}_{i_1 i_2 i_4}(\la)  \ldots D^{i_1  \ldots i_{s-1} i_s}_{i_1  \ldots i_{s-1} i_{s+1}}(\la)}\right)^2.
 \end{align}
Factorising the nominator and denominator of the fraction
on the right hand side of \eqref{A8} into the products of irreducible polynomials, rewriting \eqref{A8} in the form
\begin{align}&
\nn
\big(D^{i_1 i_2 i_3}_{i_1 i_2 i_4}(\la)  \ldots D^{i_1  \ldots i_{s-1} i_s}_{i_1  \ldots i_{s-1} i_{s+1}}(\la)\big)^2\det G(\la)
\\ \label{A9}  &
=\big(D^{i_1}_{i_2}(\la)  D^{i_1  \ldots i_{4} i_5}_{i_1  \ldots i_{4} i_{6}}(\la)\ldots 
D^{i_1  \ldots i_{r-1} i_r}_{i_1  \ldots i_{r-1} i_{r+1}}(\la)\big)^2,
 \end{align}
 and taking into account the uniqueness of the factorisation, we see that $\det G(\la)$ is a square
 of some polynomial, that is, there is a polynomial $\sqrt{\det G(\la)}$.
 Moreover, \eqref{A8} yields \eqref{A6}. 
 
 Similarly to $\det G(\la)$, one can deal with $D^{i_1   i_{2}}_{i_1   i_{2}}(\la)$
 and show that 
\begin{align}&
\label{A10}
D^{i_1   i_{2}}_{i_1   i_{2}}(\la)=p(\la)^2
\end{align} 
for some polynomial $p(\la)$. According to relations \eqref{A7} and \eqref{A10} and to the first formula in \eqref{A4},
we have
\begin{align}&
\label{A11}
D^{i_1 }_{i_1}(\la)\big/ \sqrt{\det G(\la)}=0; \quad D^{i_1 }_{i_2}(\la)\big/ \sqrt{\det G(\la)}=\breve \a p(\la),
\end{align}
where $\breve \a$ equals either $1$ or $-1$. Clearly, the expressions on the left hand 
side of the equalities in \eqref{A11}  coincide (up to the sign)
with the entries of $\sqrt{\det G(\la)}\,G(\la)^{-1}$ (and all the entries may be obtained by the
proper choice of $i_1$ and $i_2$). Therefore, the entries of $\sqrt{\det G(\la)}\,G(\la)^{-1}$ are, indeed, polynomials.
\end{proof}
\begin{La}\label{LaTh}
Let the conditions of Theorem \ref{TmMom} hold and let $\wt \th$ be given by \eqref{T80}.
Then we have
 \begin{align}\label{A12}&
\theta(\la)=-q(\la)\big(1+\wt \theta(\la)\big),
 \end{align}
where $q$ is given in \eqref{S7}. 
\end{La}
\begin{proof}. In view of \eqref{T5}, \eqref{T6} and \eqref{T34}, $q(\la)\wh u(\la)$  is a vector polynomial.
Furthermore, if  $q(\la)\wh u(\la)$ turns at some point $\la = c$ to zero, then formulas \eqref{T9}--\eqref{T6} and \eqref{T32}
yield
\begin{align}&
\label{A13}
\1^*(A_1^*-\mu_1 I)^{-1}(A_2^*-\mu_2 I)^{-1}T^{-1}h_c\equiv 0 \quad {\mathrm{for}} \quad h_c \not=0 \quad (h_c \in \BC^{mn}).
\end{align} 
However, \eqref{A13} is impossible for any nonzero $h_c$ because the span of \\
$\1^*(A_1^*-\mu_1 I)^{-1}(A_2^*-\mu_2 I)^{-1}$
coincides with $\BC^{mn}$ (and $\det T^{-1}\not=0$). Hence, $q\wh u$ does not have zeros.

It follows from \eqref{T83'} that
\begin{align}&
\label{A14}
\big(q\wh u\big)(\la)=P(\la)G(\la)^{-1}\col \begin{bmatrix}0 & \1_m &0 & \1_n\end{bmatrix}, \quad P(\la):=\I q(\la)(1+\wt \th(\la)),
\end{align} 
where (according to \eqref{T80}) $P(\la)$ is a polynomial of degree $m+n$ with the coefficient $\I$ before $\la_1^n \la_2^m$. 
Since $P(\la)G(\la)^{-1}\col \begin{bmatrix}0 & \1_m &0 & \1_n\end{bmatrix}$  is a vector polynomial and does not have zeros,
all other polynomials $\wt P(\la)$, such that $\wt P(\la)G(\la)^{-1}\col \begin{bmatrix}0 & \1_m &0 & \1_n\end{bmatrix}$
is a vector polynomial, may be factored into the product of $P$ and some other polynomial. On the other hand, it is stated in Lemma
\ref{LaSqrt} that $\sqrt{\det G(\la)}\,G(\la)^{-1}$ is a vector polynomial, $\th(\la)=\sqrt{\det G(\la)}$ is a polynomial such that
the coefficient before $\la_1^n \la_2^m$ equals $-1$ (and the degree of $\sqrt{\det G(\la)}$ equals $m+n$).
It is immediate that $\sqrt{\det G(\la)}=\I P(\la)$. Now, \eqref{A12} follows from the definition of $P$ in \eqref{A14}.

\end{proof}

\begin{flushright}

A.L. Sakhnovich,\\
Fakult\"at f\"ur Mathematik, Universit\"at Wien, \\
Oskar-Morgenstern-Platz 1, A-1090 Vienna, Austria\\
E-mail: oleksandr.sakhnovych@univie.ac.at
\end{flushright}


\begin{thebibliography}{AGKS}
\bibitem{Ald}
\textit{A. Aldroubi, J. Davis} and \textit{ I. Krishtal},  Exact reconstruction of signals in evolutionary systems via spatiotemporal trade-off,  J. Fourier Anal. Appl. textbf{21} (2015), no. 1, 11--31.

\bibitem{AGKLS}
\textit{D. Alpay, I.  Gohberg, M.A. Kaashoek,  L. Lerer} and \textit{A.L. Sakhnovich}, 
Krein systems and canonical systems on a finite interval: accelerants with a jump discontinuity at the origin and continuous potentials, 
Integral Equations Operator Theory \textbf{68} (2010), no. 1, 115--150.

\bibitem{Amb}
\textit{V.A. Ambartsumyan}, Scientific  Works (Russian). I:  Nauka, Erevan, 1960.

\bibitem{Bin}
\textit{D, Bini, S. Dendievel, G. Latouche} and  \textit{B. Meini},  General solution of the Poisson equation for quasi-birth-and-death processes, SIAM J. Appl. Math.  \textbf{76} (2016), no. 6, 2397--2417.

\bibitem{BogoB} 
\textit{J.M. Bogoya, A. B\"ottcher,  S.M.  Grudsky} and  \textit{E.A.  Maximenko},  Eigenvectors of Hermitian Toeplitz matrices with smooth simple-loop symbols, Linear Algebra Appl.  \textbf{493} (2016), 606--637.

\bibitem{BoG} 
\textit{A. B\"ottcher} and  \textit{S.M. Grudsky},  Spectral properties of banded Toeplitz matrices,  Society for Industrial and Applied Mathematics (SIAM), Philadelphia 2005.

\bibitem{BoeS} 
\textit{A.  B\"ottcher} and  \textit{B. Silbermann}, 
Analysis of Toeplitz operators, 
second edition (prepared jointly with Alexei Karlovich),  Springer
Monographs in Mathematics,  Springer-Verlag, Berlin 2006.

\bibitem{BrH} 
\textit{A. Brown} and  \textit{P.R. Halmos}, 
Algebraic properties of Toeplitz operators,  J. Reine Angew. Math. \textbf{213} (1963/1964),  89--102. 

\bibitem{CrL}
\textit{D.G. Crowdy} and  \textit{E. Luca} , Solving Wiener--Hopf problems without kernel factorization, Proc. R. Soc. Lond. Ser. A Math. Phys. Eng. Sci. \textbf{470} (2014), no. 2170, 20140304, 20 pp.

\bibitem{FKRS}
\textit{B. Fritzsche, B. Kirstein, I.Ya. Roitberg} and  \textit{A.L. Sakhnovich},  Operator identities corresponding to inverse problems, Indag. Math., New Ser.  \textbf{23} (2012), no. 4, 690--700. 

\bibitem{GaKr}
\textit{F.P.  Gantmacher} and  \textit{M.G. Krein},  Oscillation matrices and kernels and small vibrations of mechanical systems, GITTL, Moscow 1950.  Translation edited by A. Eremenko, AMS Chelsea Publishing, Providence 2002.

\bibitem{GeK}
\textit{Y. Genin} and  \textit{Y. Kamp}, Two-dimensional stability and orthogonal polynomials on the hypercircle, in:  Proc. IEEE \textbf{65} (1977), no. 6, 873--881. 


\bibitem{Go}
\textit{I. Gohberg}  (Ed.), {Continuous and discrete Fourier transforms, extension problems and Wiener--Hopf equations}, Operator Theory Adv. Appl. \textbf{58}, Birkh\"auser, Basel 1992.  

\bibitem{GoKr}
\textit{I.C. Gohberg} and  \textit{M.G. Krein}, Systems of integral equations on a half line with kernels depending on the difference of arguments, Amer. Math. Soc. Transl. (2) \textbf{14} (1960), 217--287.

\bibitem{GoS} 
\textit{I.C. Gohberg} and  \textit{A.A. Semencul},   The inversion of finite Toeplitz matrices and their continual analogues (Russian),  Mat. Issled. \textbf{7} (1972), no. 2, 201--223.

\bibitem{Gr} 
\textit{K. Gr\"ochenig} and  \textit{M. Leinert}, 
Wiener's lemma for twisted convolution and Gabor frames, 
J. Amer. Math. Soc. \textbf{17} (2004), no. 1, 1--18.

\bibitem{JeVa}
\textit{B. Jeuris} and  \textit{R. Vandebril},  The K\"ahler mean of block-Toeplitz matrices with Toeplitz structured blocks, SIAM J. Matrix Anal. Appl. \textbf{37} (2016), no. 3, 1151--1175.

\bibitem{Jus}
\textit{J.H. Justice}, A Levinson type algorithm for two dimensional Wiener fittering using bivariate Szeg\"o polynomials, in: Proc. IEEE \textbf{65} (1977), no. 6, 882--886.

\bibitem{Kalo}
\textit{N. Kalouptsidis, G. Carayannis} and  \textit{D. Manolakis},  On block matrices with elements of special structure, in: Proceedings of IEEE International Conference on Acoustics, Speech and Signal Processing (ICASSP-82) \textbf{7},  IEEE, New York 1982, 1744--1747.

\bibitem{Kre}
\textit{M.G. Krein},  Integral equations on the half-line with a kernel depending on the difference of the arguments, Amer. Math. Soc. Transl. \textbf{22}  (1962), 163--288.  

\bibitem{Levin}
\textit{N. Levinson}, The Wiener rms (root mean square) error criterion
in filter design and prediction, J. Math. Phys. \textbf{25} (1947),  261--278.

\bibitem{Lev}
\textit{B. Levy}, $2-D$ polynomial and rational matrices and their application for the modeling of $2-D$ dynamical systems, Ph.D. dissertation, Stanford Univ., Stanford, California 1981.

\bibitem{Nap}
\textit{V.V.  Napalkov},  Convolution equations in multidimensional spaces (Russian), Nauka, Moscow 1982.

\bibitem{NiB} 
\textit{F. Nielsen} and  \textit{R. Bhatia} (Eds), Matrix information geometry, Springer-Verlag, Heidelberg 2013. 

\bibitem{Nob} 
\textit{B. Noble},  Methods based on the Wiener--Hopf technique,  Pergamon Press, London 1958.

\bibitem{SaA73} 
\textit{A.L. Sakhnovich}, On a method of inverting Toeplitz  matrices, Math. Issled. \textbf{8} (1973), no. 4, 180--186.

\bibitem{SaAJFA} 
\textit{A.L. Sakhnovich},  Toeplitz matrices with an exponential growth of entries and the first Szeg\"o limit theorem, J. Functional Anal. \textbf{171} (2000), 449--482.

\bibitem{SaA2017} 
\textit{A.L. Sakhnovich}, Inversion of the convolution operators on a rectangular,
 arXiv:1701.08559

\bibitem{SaSaR}
\textit{A.L.  Sakhnovich, L.A. Sakhnovich} and  \textit{I.Ya. Roitberg}, Inverse problems and nonlinear evolution equations. 
 Solutions, Darboux matrices and Weyl--Titchmarsh functions, De Gruyter Studies in Mathematics \textbf{47}, De Gruyter, Berlin 2013.  

\bibitem{SaL68}
\textit{L.A. Sakhnovich},
 Operators, similar to unitary operators, with absolutely continuous spectrum (Russian),  Funkcional. Anal. i Prilozhen.  \textbf{2} (1968), no.1, 51--63. 

\bibitem{SaL73}
\textit{L.A. Sakhnovich},
An integral equation with a kernel dependent on the difference of the arguments, Mat. Issled.  \textbf{8} (1973), no. 2, 138--146.


 \bibitem{SaL80}
\textit{L.A. Sakhnovich}, Equations with a difference kernel on a  finite interval,
Russ. Math. Surveys \textbf{35} (1980), no. 4,   81--152.

\bibitem{SaL15}
\textit{L.A. Sakhnovich},  \textit{Integral equations with difference kernels on finite intervals}, second edition, revised and extended,
 Operator Theory Adv. Appl. \textbf{84}, Birkh\"auser/Springer, Cham 2015.

 \bibitem{Sha} 
\textit{D. Sharpe}, Rings and factorization,  Cambridge University Press, Cambridge 1987.

\bibitem{Shi}
\textit{M. Shinbrot}, A class of difference kernels, Proc. Amer. Math. Soc. \textbf{13} (1962), 399--406. 

 \bibitem{Tr}
\textit{W.F. Trench}, An algorithm for the inversion of prediction of stationary time series from finite past, SIAM J. Appl. Math. \textbf{15} (1967),  no. 6, 1502--1510.

\bibitem{WaK} 
\textit{M. Wax} and  \textit{T. Kailath}, Efficient inversion of Toeplitz-block Toeplitz matrix,
 IEEE Trans. Acoust. Speech Signal Process. \textbf{31} (1983), no. 5, 1218--1221.


\bibitem{Xi}
\textit{P. Xie} and  \textit{Y. Wei},  The stability of formulae of the Gohberg--Semencul--Trench type for 
Moore--Penrose and group inverses of Toeplitz matrices,  Linear Algebra Appl. \textbf{498} (2016), 117--135.


\end{thebibliography}
\end{document}